# Quantum Optimisation for Transport Vulnerability Identification

Junxiang Xu[1,*], Chence Niu[2], Divya Jayakumar Nair[1], Vinayak Dixit[1]

*1. Research Centre for Integrated Transport Innovation (rCITI), School of Civil and Environmental Engineering, The University of New South Wales, Kensington, UNSW Sydney, NSW, 2052, Australia*

*2. Guangdong Basic Research Center of Excellence for Ecological Security and Green Development, Key Laboratory for City Cluster Environmental Safety and Green Development of the Ministry of Education, School of Ecology, Environment and Resources, Guangdong University of Technology, Guangzhou, 510006, China*

**corresponding author*:** junxiang.xu@unsw.edu.au

**Abstract:** Transport network vulnerability analysis plays a crucial role in safeguarding urban resilience. Traditional vulnerability identification approaches have provided valuable insights, yet they face two major limitations. First, the number of disruption scenarios increases combinatorially with the number of disrupted links considered simultaneously, making classical approaches computationally prohibitive. Second, most studies approximate the impacts of multiple simultaneous link failures through linear aggregation, which fails to capture the nonlinear interaction effects observed in real networks. To address these gaps, we reformulate the bi-level Mixed-Integer Nonlinear Programming (MINLP) model into a quantum-compatible Quadratic Unconstrained Binary Optimisation (QUBO) structure, enabling parallel exploration of complex disruption scenarios while incorporating nonlinear interaction effects. We develop a hybrid optimisation framework that integrates the quantum optimisation algorithm with the Frank-Wolfe method to validate the model's effectiveness on the small-scale network. Then, we further verify the framework through the D-Wave hardware across benchmark networks of different scales, including the Sioux Falls, Anaheim, Chicago Sketch, and Berlin Full networks, to examine its scalability and hardware feasibility. The results show that this framework achieves strong solvability and stability. In particular, optimisation for large and larger networks is completed within minutes (Approximately 2.8 minutes for the 914-link, 9.8 minutes for the 2950-link, and 31.2 minutes for the 6018-link on D-Wave), demonstrating a computational efficiency improvement by one to two orders of magnitude compared with classical metaheuristic algorithms. These findings highlight the feasibility and potential of applying quantum computing to network vulnerability identification and open a new avenue for resilience-oriented planning.

## 1. Introduction

With accelerating urbanisation and the continuous growth of travel demand, transport networks play a vital role in modern society. However, they are also highly vulnerable, as the failure of certain critical links due to natural hazards, human-induced events or infrastructure deterioration can trigger cascading effects and lead to severe degradation of overall network performance (Rodríguez-Núñez and García-Palomares, 2014; Chen et al., 2015; Mattsson and Jenelius, 2015; Bell et al., 2017; Mouronte-López, 2021). Numerous studies have shown that the disruption of only a small number of critical links can paralyse the efficiency of the entire network, creating risks for economic activity, environmental sustainability, and social stability (Scott et al., 2006; Sullivan et al., 2010; Wang et al., 2016; Zhou and Wang, 2017; Jafino et al., 2020; Du et al., 2022). Consequently, how to effectively identify and evaluate the vulnerability of transport networks, particularly the vulnerability of critical links, has become a central issue in the fields of transport research and urban planning.

In recent years, the rapid development of quantum computing has opened new possibilities for addressing complex combinatorial optimisation problems (Khumalo et al., 2022; Heng et al., 2022; Chicano et al., 2025). Unlike traditional computational methods, which rely on exhaustive enumeration or heuristic search, quantum computing exploits unique features such as superposition, entanglement and quantum tunnelling to explore vast solution spaces in parallel, and offers the potential to overcome the scalability bottlenecks of conventional algorithms (Khurana et al., 2024; Zhao et al., 2024; Jarju, 2025). In particular, the combination of quantum annealing with Quadratic Unconstrained Binary Optimisation (QUBO) modelling enables a class of inherently NP-hard problems to be reformulated within a quantum framework for more efficient solution. Existing studies have already demonstrated the promise of quantum algorithms in fields such as financial portfolio optimisation (Herman et al., 2023), energy scheduling (Mastroianni et al., 2023), logistics and complex network design (Niu et al., 2025b). Quantum computing holds significant potential for large-scale applications in transport network vulnerability identification and resilience enhancement.

Therefore, this study integrates quantum computing with the problem of vulnerability identification in transport resilience, exploring a new pathway to overcome the limitations of classical approaches and address the scalability challenges arising from the combinatorial explosion of disruption scenarios. By reformulating the classical bi-level vulnerability identification model into a quantum-compatible QUBO structure and designing a hybrid quantum-classical algorithmic framework, this study provides a novel modelling and solution approach for identifying critical links in transport networks. The main contributions are threefold:

(1) This study is the first to reformulate the bi-level optimisation model for critical link identification into a quantum-compatible QUBO form, enabling integration with quantum computing and capturing nonlinear interaction effects among multiple link disruptions.

(2) This study develops a hybrid quantum-classical framework combining Simulated Quantum Annealing (SQA), Quantum Annealing (QA), and the Frank-Wolfe (FW) algorithm, ensuring solvability and stability while enhancing exploration efficiency.

(3) This study further implements the framework on real quantum hardware using the D-Wave platform and validates scalability across the Nguyen-Dupuis, Sioux Falls, Anaheim, Chicago Sketch, and Berlin Full networks, demonstrating that large-scale optimisation can be completed within minutes, far outperforming traditional metaheuristic algorithms.

The rest of this paper is arranged as follows. **Section 2** is the literature review. **Section 3** presents the classical bi-level formulation of the critical link identification problem. **Section 4** introduces the quantum-compatible reformulation and parameter specification. **Section 5** designs a hybrid quantum-classical optimisation algorithm. **Section 6** first presents experiment on the ND network and then extends the validation to medium and large-scale networks on D-Wave quantum hardware. **Section 7** concludes this paper.

## 2. Literature review

Transport network vulnerability refers to the sensitivity of a transport system to disturbances and the extent of performance degradation, and it has become one of the core research topics within the broader study of transport resilience (Mattsson and Jenelius, 2015; Reggiani et al., 2015). Transport network vulnerability is widely recognised as a key dimension of resilience analysis. Together with robustness and connectivity, it forms part of the broader framework of network resilience and reliability research (Xu and Nair, 2024). According to existing studies, research on transport network vulnerability can be categorised into three main areas. First, vulnerability measurement, which focuses on quantifying performance degradation when certain links or nodes fail, using indicators such as network connectivity, reliability indices, travel time, accessibility and travel cost (Chen et al., 2007; Du et al., 2014; Lu et al., 2021). Second, critical link or node identification, which examines the extent to which the failure of a small number of links or nodes impacts the overall network performance, this line of research often directly informs resilience assessment and the prioritisation of resilience investments (Kumar et al., 2019; Guettiche and Kheddouci, 2019; Lai and Zhang, 2022). Third, multi-scenario disruption and multi-objective optimisation, which develops multi-objective optimisation models under conditions such as simultaneous link failures, demand uncertainty or different disaster scenarios, to provide a more comprehensive evaluation of transport network vulnerability and resilience (Acosta et al., 2018; Pérez-Morales et al., 2019; Tiong and Vergara, 2023; Jiang and Song, 2024).

In terms of transport network vulnerability research approaches, existing studies can be categorised into three. The first type comprises methods based on network topology or complex network structural attributes, such as metrics including betweenness centrality, connectivity loss and the number of alternative paths (Reggiani et al., 2015; Zhang et al., 2015; Chen and Lu, 2020; Chen et al., 2023). These methods are computationally efficient, do not rely on detailed network attribute data such as traffic flow, travel time or cost, but they often fail to capture the actual characteristics for traffic flows. The second includes bi-level optimisation models in transport network design theory. In such models, the upper level is formulated as a vulnerability optimisation problem by introducing binary decision variables to identify whether a road link or node is critical, while the lower level corresponds to a User Equilibrium (UE) or Stochastic User Equilibrium (SUE) (Zhang et al., 2020; Gan et al., 2021; Ghanaei et al., 2025). This approach can more realistically represent the impact of road link failures on overall system performance, though it involves relatively complex modelling and solution procedures. The third method refers to simulation and statistical methods, which use Monte Carlo simulation or random scenario generation techniques to calculate the distribution of network performance under multiple disruption cases, thereby assessing the level of vulnerability and potential risks (Praks et al., 2017; Haghighi et al., 2018; Chen et al., 2021; Liu et al., 2024). There is a common challenge across all these methods is that, as the number of disrupted links increases, the combinatorial complexity grows exponentially, resulting in extremely high computational costs and limiting the applicability of such models to large-scale real-world transport networks.

Within this broader stream of research, critical link identification has emerged as a central topic, as the failure of a small number of links often leads to disproportionate increases in total system travel time and severe degradation of network performance. Various approaches have been proposed to address this problem, including bi-level optimisation frameworks for identifying vulnerable link combinations (Li et al., 2017), upper and lower bound formulations for performance under multiple disruptions (Xu et al., 2017), user equilibrium interdiction models (Starita and Scaparra, 2021), analyses of critical link combinations in urban road networks (Jin et al., 2022), and ranking methodologies based on criticality score distributions (Barati et al., 2024). These studies have substantially enriched the theoretical and methodological foundations of vulnerability analysis and consistently demonstrated that network fragility is highly concentrated in a small set of critical links. Nevertheless, they remain constrained by the computational bottlenecks discussed above, as the combinatorial explosion associated with multiple-link disruptions severely limits their applicability to large-scale real-world transport networks.

As highlighted in the Introduction section, quantum computing has in recent years demonstrated potential in areas such as financial optimisation, energy scheduling and supply chain management. However, its application in the transport domain has developed more slowly, with existing studies focusing mainly on transport network design, resilience-oriented investment, travelling salesman problem, and multi-commodity flow problems (Dixit and Jian, 2022; Dixit and Niu, 2023; Dixit et al., 2024; Ali et al., 2024; Jung and Choi, 2025; Smith-Miles et al., 2025; Sato et al., 2025; Niu et al., 2025b). Most of these studies reformulate complex transport optimisation models into QUBO structures and employ quantum annealing, the Quantum Approximate Optimisation Algorithm (QAOA) or hybrid quantum-classical frameworks to obtain solutions (Salloum et al., 2025). At the current stage, algorithmic applications of quantum computing in transport generally include gate-based models, quantum annealing and quantum machine learning (Niu et al., 2025a). Overall, the principal role of quantum algorithms in transport research lies in enhancing computational efficiency through parallel search and quantum tunnelling, thereby alleviating, to some extent, the combinatorial explosion encountered by classical methods in networks of certain scales (Cooper, 2021; Gabbassov, 2022; Dixit and Niu, 2023). Although research in this area is still at an early stage, these developments open new possibilities for applying quantum computing to transport network vulnerability identification and provide a promising direction for resilience-oriented transport planning.

In summary, research on transport network vulnerability identification has accumulated substantial achievements in vulnerability measurement, critical link identification and multi-scenario disruption analysis. However, existing approaches largely rely on mathematical programming, exact algorithms or evolutionary heuristics, and are typically designed to identify critical links under random or deliberate disruptions. In real-world networks, multiple links may fail simultaneously, and most existing studies approximate the consequences of such failures through linear aggregation, which limits their ability to capture potential nonlinear interaction effects among disrupted links (Iliopoulou and Makridis, 2023; Chen et al., 2023). As a result, when considering disruption scenarios involving a certain number of simultaneous link failures, current methods face significant challenges in both modelling and computation (Wang et al., 2024; Boeing and Ha, 2024). At the same time, the potential of quantum computing in large-scale combinatorial optimisation has already been demonstrated, and applications have begun to emerge in related areas such as transport network design and resilience-oriented investment. Features such as quantum superposition and quantum tunnelling provide new modelling opportunities to more realistically represent the complex impacts of multiple link disruptions (Wen et al., 2024). Nevertheless, there remains a lack of systematic efforts to introduce quantum modelling and quantum algorithms into the transport domain, particularly in relation to vulnerability identification. Addressing this gap, this study seeks to integrate quantum computing with the problem of critical link identification in transport networks, reformulating combinatorially complex optimisation models into quantum-compatible QUBO structures and developing a hybrid quantum-classical framework that provides methodological support for resilience-oriented transport planning.

Therefore, this study aims to address three core questions:

*RQ1: How can the classical transport network vulnerability identification model be reformulated into a quantum-compatible QUBO structure suitable for quantum optimisation?*

*RQ2: How can quantum computing capture and evaluate the nonlinear interaction effects of multiple simultaneous link disruptions across different network scales?*

*RQ3: How can a hybrid quantum-classical algorithm be designed, implemented, and validated on real quantum hardware to ensure solvability, convergence, and scalability in large-scale networks?*

The technical approach for this paper is designed in **Appendix A**.

### 3. Methodology

To facilitate the introduction of the methodology in this section, especially the modelling details, the parameter settings shown in **Table 1** are defined first.

**Table 1.** Parameter settings and definitions.

| Parameter | Description |
|---|---|
| $S$ | Set of road links in the network, $S=\{1,2,...,s\}$ |
| $W$ | Set of origin-destination (OD) pairs, $W=\{1,2,...,\omega\}$ |
| $R$ | Set of all feasible routes in the network, $R=\{1,2,...,r\}$ |
| $d^{\omega}$ | Travel demand of OD pair $\omega$ |
| $C_s$ | Original capacity of road link $s$ |
| $t_s^0$ | Free-flow travel time of road link $s$ |
| $\alpha$ | BPR function parameter (commonly set as 0.15) |
| $\beta$ | BPR function parameter (commonly set as 4.0) |
| $e$ | Remaining capacity ratio after disruption, $e\in[0,1]$ |
| $k$ | Number of simultaneously disrupted road links |
| $u_s$ | Binary decision variable: $u_s=1$ if road link $s$ is disrupted, otherwise 0 |
| $x_s$ | Traffic flow on road link $s$ |
| $f_r^{\omega}$ | Traffic flow of OD pair $\omega$ on route $r$ |
| $\tau_{sr}$ | Route-link incidence variable: equals 1 if route $r$ contains road link $s$ , otherwise 0 |
| $T_s(x_s,u_s)$ | Travel time function of road link $s$ (extended BPR function) |

This study formulates the critical road identification problem in transport networks as a bi-level Mixed-Integer Nonlinear Programming (MINLP) model (Li et al., 2017; Xu et al., 2017; Starita and Scaparra, 2021; Jin et al., 2022; Barati et al., 2024). In this framework, the upper level represents the disruption strategy, i.e., selecting a fixed number of road links to be disrupted to maximise the Total System Travel Time (TSTT). The lower level corresponds to the UE traffic assignment problem, where road users respond to the given disruptions by redistributing their trips across the network to minimise their individual travel costs.

For the sake of tractable modelling, several assumptions are made in this study:

(1) The travel demand between each OD pair is assumed to be fixed and inelastic, regardless of road disruptions.

(2) Travellers are rational and follow Wardrop's first principle, such that the traffic assignment reaches UE (Wardrop and Whitehead, 1952).

(3) The travel time of each road link follows the extended Bureau of Public Roads (BPR) function, where the capacity reduction due to disruptions is represented by a fixed proportion of the original capacity. (Basic BPR function: $T_s(x_s) = t_s^0\left[1+\alpha\left(\frac{x_s}{C_s}\right)^{\beta}\right], \forall s \in S$ (Weiner, 1987)

(4) The disruption is complete and static, where a road is either in a normal state or its capacity is reduced to $eu_s$ , without considering partial repair or dynamic recovery.

(5) It is assumed that the road network remains connected even under disruptions, i.e., at least one feasible route exists between every OD pair.

(6) Stochastic variations in traffic demand, travel behaviour, or random incidents are not considered. All parameters are treated as deterministic.

Based on the above assumptions and parameter settings, the classical bi-level MINLP optimisation model can be established as follows.

Upper-level model (critical road disruption identification):

$$Z = \max \sum_{s \in S} x_s T(x_s, u_s) \tag{1}$$

subject to:

$$\sum_{s \in S} u_s = k \tag{2}$$

$$T_s\left(x_s,u_s\right)=t_s^0\left[1+\alpha\left(\frac{x_s}{C_s\left(1-u_s+eu_s\right)}\right)^{\beta}\right],\forall s\in S \tag{3}$$

$$u_s\in\{0,1\},\forall s\in S \tag{4}$$

Lower-level model (user equilibrium traffic assignment):

$$F=\min\sum_{s\in S}\int_0^{x_s}T_s\left(\varphi,u_s\right)d\varphi \tag{5}$$

subject to:

$$d^{\omega}=\sum_{r\in R}f_r^{\omega},\forall\omega\in W \tag{6}$$

$$x_s=\sum_{\omega\in W}\sum_{r\in R}f_r^{\omega}\tau_{sr},\forall s\in S \tag{7}$$

$$x_s\geq 0,f_r^{\omega}\geq 0,\tau_{sr}\in\{0,1\},\forall s\in S,r\in R,\omega\in W \tag{8}$$

**Eq. (1)** is the upper-level objective function, which maximises the TSTT under road disruption scenarios, to characterise the most severe degradation of network performance. **Eq. (2)** constrains the number of simultaneously disrupted road links to be exactly $k$. **Eq. (3)** defines the road travel time function (an extended BPR function), which incorporates the effect of capacity reduction under road disruption scenarios. **Eq. (4)** is the binary decision variable, indicating whether a road is disrupted. **Eq. (5)** is the lower-level objective function, namely the UE traffic assignment model, which characterises the equilibrium flow distribution under Wardrop's first principle by minimising the total travel cost of the network. **Eq. (6)** is the OD demand conservation constraint, ensuring that the demand of each OD pair is fully assigned to the network routes. **Eq. (7)** is the relationship between road flows and route flows, ensuring that the flow on each road equals the sum of the flows of all routes traversing that road. **Eq. (8)** represents the non-negativity and binary constraints.

## 4. Quantum modelling

### *4.1 Reformulation of the bilevel MINLP*

Classical algorithms for bi-level MINLP problems, such as branch-and-bound, cutting-plane, and heuristic search, are well developed (Lachhwani and Dwivedi, 2018; Caselli et al., 2024; Liñán and Ricardez-Sandoval, 2025; Martens, 2025; Hammad et al., 2025). They can solve disruption identification problems in small and medium-scale networks. However, their computational complexity grows rapidly with network size and the number of disrupted road links. Quantum computing provides an alternative. Techniques such as quantum annealing, quantum machine learning, gate-based quantum computing, and variational quantum algorithms can explore large solution spaces in parallel (Niu et al., 2025a). They are also more effective in avoiding local optimisation. In this study, the upper-level disruption identification problem is reformulated into a QUBO structure for quantum solvers. The lower-level UE model remains in a classical framework to preserve established traffic flow principles. This hybrid quantum-classical design combines the scalability of quantum optimisation with the stability of classical models.

To transfer the upper-level disruption identification model into a quantum-compatible model, the original MINLP expression needs to undergo binarization and quadratic transformation. The original objective and constraints are:

$$Z = \max \sum_{s \in S} x_s T(x_s, u_s), \sum_{s \in S} u_s = k, u_s \in \{0,1\}, \forall s \in S \tag{9}$$

As can be seen in **Table 1**, $u_s$ denotes the binary decision variable indicating whether road link $s$ is disrupted. Since qubits are inherently binary with values $\{0,1\}$, $u_s$ can be directly mapped to qubits. Meanwhile, the constraint $\sum_{s \in S} u_s = k$ can be incorporated into the objective function as a quadratic penalty term.

To adapt the model formulations for quantum optimisation algorithms, which generally operate on QUBO structures, the objective is reformulated as a minimisation problem in QUBO form. Following the general QUBO formulation principle in quantum optimisation (Lucas, 2014; Glover et al., 2018; Glover et al., 2022), the QUBO formulation in this study is expressed as follows.

$$Z_{\text{qubo}} = \min H(u) = -\sum_{s \in S} c_s u_s - \sum_{s,t \in S} \beta_{st} u_s u_t + \lambda \left( \sum_{s \in S} u_s - k \right)^2 \tag{10}$$

where $c_s$ represents the impact of disrupting road link $s$ on the total system travel time (TSTT), derived from the lower-level UE assignment model. $\beta_{st}$ is an interaction coefficient between road links, introduced to enhance quadratic representation efficiency. $\lambda$ is a penalty coefficient ensuring that the constraint $\sum_{s \in S} u_s = k$ is satisfied.

Then, rewrite the model as follows (Hybrid Quantum-Classical Bi-level Model)
Upper-level model (critical road disruption identification, QUBO form):

$$Z_{\text{qubo}} = \min H(u) = -\sum_{s \in S} c_s u_s - \sum_{s,t \in S} \beta_{st} u_s u_t + \lambda \left( \sum_{s \in S} u_s - k \right)^2 \tag{10}$$

Lower-level model (user equilibrium traffic assignment):

$$F = \min \sum_{s \in S} \int_0^{x_s} T_s(\varphi, u_s) \, d\varphi \tag{5}$$

subject to:

$$d^{\omega} = \sum_{r \in R} f_r^{\omega}, \forall \omega \in W \tag{6}$$

$$x_s = \sum_{\omega \in W} \sum_{r \in R} f_r^{\omega} \tau_{sr}, \forall s \in S \tag{7}$$

$$x_s \geq 0, f_r^{\omega} \geq 0, \tau_{sr} \in \{0,1\}, \forall s \in S, r \in R, \omega \in W \tag{8}$$

*4.2 Parameter specification in the QUBO formulation*

Parameter specification in the QUBO formulation:

(1) Single-link impact coefficient $c_s$

Specifically, $c_s$ measures the marginal increase in the TSTT when the road link $s$ is disrupted. Firstly, solve the lower-level UE model twice, with $u_s = 1$ (disrupted) and $u_s = 0$ (normal). Secondly, obtain the corresponding system travel times $TSTT(u_s = 1)$ and $TSTT(u_s = 0)$. Finally, define:

$$c_s = TSTT(u_s = 1) - TSTT(u_s = 0) \tag{11}$$

$$TSTT = \sum_{s \in S} x_s T(x_s, u_s) = \sum_{s \in S} x_s t_s^0 \left[ 1 + \alpha \left( \frac{x_s}{C_s(1 - u_s + eu_s)} \right)^{\beta} \right], \forall s \in S \tag{12}$$

Eq. (11) and Eq. (12) represent the additional delay caused by disrupting road link $s$.

(2) Interaction coefficient $\beta_{st}$

Specifically, $\beta_{st}$ captures the nonlinear effect when two road links $s$ and $t$ are disrupted simultaneously, beyond the additive effects of individual disruptions. Firstly, run the UE model for the joint disruption case $u_s = 1, u_t = 1$. Then, calculate the joint TSTT and compare with the sum of single disruptions:

$$\beta_{st} = TSTT(u_s = 1, u_t = 1) - \left[ TSTT(u_s = 1) + TSTT(u_t = 1) - TSTT(u_s = 0, u_t = 0) \right] \tag{13}$$

If $\beta_{st} > 0$: This indicates a negative synergistic effect, where the two failures exacerbate congestion.

If $\beta_{st} < 0$: This suggests a degree of substitutability, meaning that the impact of one disruption partly offsets the other.

If $\beta_{st} = 0$: The combined effect is exactly equal to the sum of the two independent effects, implying a purely additive relationship without interaction.

(3) Penalty coefficient $\lambda$

$\lambda$ must be sufficiently large so that violating the constraint leads to a penalty dominating the objective. A common practice is to set $\lambda$ larger than the maximum possible sum of $|c_s|$ and $|\beta_{st}|$:

$$\lambda \gg \max\left( \sum_{s \in S} |c_s|, \sum_{s,t \in S} |\beta_{st}| \right) \tag{14}$$

When performing actual calculations, to avoid infinite values, take 10 to 100 times the maximum coefficient as the value of $\lambda$.

*4.3 Equivalence and solvability of the QUBO formulation*

The reformulated QUBO-based upper-level model remains equivalent to the original bilevel critical disruption identification problem. Specifically, the three coefficients ensure that the objective function in the QUBO form faithfully represents the marginal and joint impacts of road disruptions together with the cardinality constraint. Hence, solving the QUBO problem is mathematically consistent with solving the original combinatorial optimisation model.

From a computational perspective, the QUBO formulation transforms the bilevel problem into a quadratic unconstrained binary optimisation problem, which belongs to the class of NP-hard problems. However, this form is directly compatible with quantum optimisation paradigms such as quantum annealing and the Quantum Approximate Optimisation Algorithm (QAOA), as well as classical heuristics designed for QUBO (Dixit and Niu, 2023; Niu et al., 2025b). This ensures that the proposed hybrid quantum-classical bilevel framework is not only theoretically valid but also practically solvable using emerging quantum technologies in combination with conventional algorithms.

## 5. Algorithm design: SQA-FW and its extension to D-Wave quantum hardware

### *5.1 Upper-level QUBO optimisation using quantum annealing*

To solve the upper-level QUBO optimisation problem:

$$Z_{\text{qubo}} = \min H(u) = -\sum_{s\in S} c_s u_s - \sum_{s,t\in S} \beta_{st} u_s u_t + \lambda \left( \sum_{s\in S} u_s - k \right)^2 \tag{10}$$

this study adopts Simulated Quantum Annealing (SQA). SQA extends classical Simulated Annealing (SA) by incorporating *imaginary-time slices* and a *quantum tunnelling* term, which enables the algorithm to escape from local minima more effectively and better approximate quantum annealing behaviour.

**Algorithm 1. SQA for QUBO**

**Input:** coefficient $c_s$, interaction coefficient $\beta_{st}$, penalty coefficient $\lambda$, initial temperature $T_0$, cooling rate $v, v \in [0,1]$, number of imaginary-time slices $M$, maximum number of iterations $N_{iter}$.

**Output:** optimal solution vector $u_s^*$ (optimised result for $u_s$ in the $Z_{\text{qubo}}$)

*Step 1. Initialisation*

Randomly generate an initial solution $u_s^{(0)} \in \{0,1\}^{|S|}$ and replicate it across $M$ imaginary-time slices. Set the initial temperature $T = T_0$ and the quantum tunnelling strength $\Gamma > 0$.

*Step 2. Energy evaluation*

(1) Calculate the QUBO energy

$$E(u) = -\sum_{s\in S} c_s u_s - \sum_{s,t\in S} \beta_{st} u_s u_t + \lambda \left( \sum_{s\in S} u_s - k \right)^2 \tag{15}$$

(2) Extend the energy across the time slices

$$E_{SQA}\left(u_s^*\right) = \sum_{m=1}^{M} E\left(u_s^{(m)}\right) - \Gamma \sum_{m=1}^{M} \sum_{s\in S} u_s^{(m)} u_s^{(m+1)} \tag{16}$$

*Step 3. State update*

(1) At each iteration, randomly select a bit $u_s^{(m)}$ and attempt to flip it.

(2) Calculate the energy difference $\Delta E_{SQA}\left(u_s^*\right)$ (Eq. (16)), and accept the new state with probability:

$$P = \min\left\{1, \exp\left(-\frac{\Delta E_{SQA}\left(u_s^*\right)}{T}\right)\right\} \tag{17}$$

*Step 4. Annealing schedule*

(1) After each iteration, update the temperature $T_{n+1} = vT_n, n = \{1,2,...,N_{iter}\}$ , The value of $v$ is usually taken as 0.95-0.99.

(2) Gradually reduce the tunnelling strength $\Gamma$ until it approaches zero.

Generally, a larger $\Gamma$ is set initially to ensure that the system has sufficient 'quantum fluctuations' to escape local optima, such as $\Gamma_0 = 10.0$ . Like temperature annealing ( $T$ from high to low), $\Gamma$ should also be gradually reduced until it approaches 0. This study uses exponential decay:

$$\Gamma(n) = \Gamma_0 \cdot e^{-vn} \tag{18}$$

Generally, $\Gamma_0$ typically chosen between 1 and 10 and adjusted in conjunction with $M$ and annealing iterations $N_{iter}$ .

This study uses the same value of $v$ to simultaneously control the synchronous speed of two annealing processes ( $T$ and $\Gamma$ ).

*Step 5. Termination and output*

Stop when the maximum iteration number $N_{iter}$ is reached or the energy stabilises.

Output the solution with the lowest energy as $u_s^*$ .

In fact, the SQA adopted in this study employs a stochastic generation method like the classical Monte Carlo approach (Kuhl, 1996; Chen et al., 1999; Chen et al., 2011), and it can be directly implemented in a Python environment. Specifically, the solution process of SQA relies on the iterative optimisation of an energy function: a feasible solution is first randomly initialised, and at each step, a bit is randomly selected for flipping. The energy difference before and after the flip is then calculated, and a modified Metropolis acceptance criterion (like simulated annealing) is applied to determine whether the new solution is accepted. As the iterations proceed, the temperature is gradually reduced according to a predefined cooling schedule, enabling the algorithm to escape local optima in the early stages and progressively converge to the global optimum in later stages. Since the entire process involves only bit operations on the solution vector, numerical evaluation of the energy function, and random number generation, it can be fully implemented using Python's numerical libraries such as NumPy. In addition, more mature simulated quantum annealers, such as the *openjij and* D-Wave, can be employed to ensure efficiency and stability.

*5.2 Lower-level traffic assignment based on user equilibrium*

The lower-level UE traffic assignment problem in Eq. (5)-(8) is solved using the Frank-Wolfe (FW) algorithm, which is a classical method for convex optimisation in transport networks (LeBlanc and Farhangian, 1981; Fukushima, 1984; Lee and Nie, 2001; Potuzak and Kolovsky, 2022; Bomze et al., 2024). The procedure is summarised as follows.

**Algorithm 2. FW for UE model**

Since the FW algorithm is so well-known, this section will not go into detail. The algorithm content and process are consistent with another existing study (Xu et al., 2024).

This procedure yields the equilibrium road link flows $x_s$, route flows $f_r^{\omega}$, and $TSTT$, which are then fed back into the upper-level model.

Based on the algorithms described in Sections 5.1 and 5.2, this section provides a detailed description of the steps involved in the hybrid optimisation algorithm, as shown in **Appendix B**.

*5.3 Extension to D-Wave hardware implementation*

Although the proposed SQA-FW algorithm employs SQA to approximate quantum effects within a classical computing environment, the upper-level optimisation can be seamlessly transferred to real quantum hardware. In this extended implementation, the simulated quantum annealing process is replaced by the D-Wave using Quantum Annealing (QA), which executes the equivalent QUBO formulation directly on quantum processing units. The FW algorithm remains unchanged at the lower level, solving the user equilibrium flows on a classical computer.

This hardware-level implementation maintains the same iterative structure of the SQA-FW algorithm while leveraging quantum hardware to explore a larger combinatorial space with potential improvements in convergence and scalability. For benchmark networks (e.g., Sioux Falls), the D-Wave Leap Hybrid Solver is adopted to combine QPU-based sampling with classical subproblem decomposition, ensuring computational feasibility while preserving quantum advantages.

## 6. Numerical experiments

*6.1 Data preparation: Nguyen-Dupuis test network*

This study adopts the Nguyen-Dupuis (ND) test network as the experimental benchmark, which has been widely used in the literature for testing traffic assignment models, network design problems and vulnerability analysis due to its moderate size and well-documented structure (Cheng et al., 2024; Tu et al., 2024; Chen et al., 2024; Nitheesh and Bhavathrathan, 2025; Xu et al., 2025). The structure of the network is illustrated in **Figure 1**, while the associated network attributes and link parameters, the OD demand matrix, and the parameters used for the proposed model and algorithm are listed in **Appendix C**.

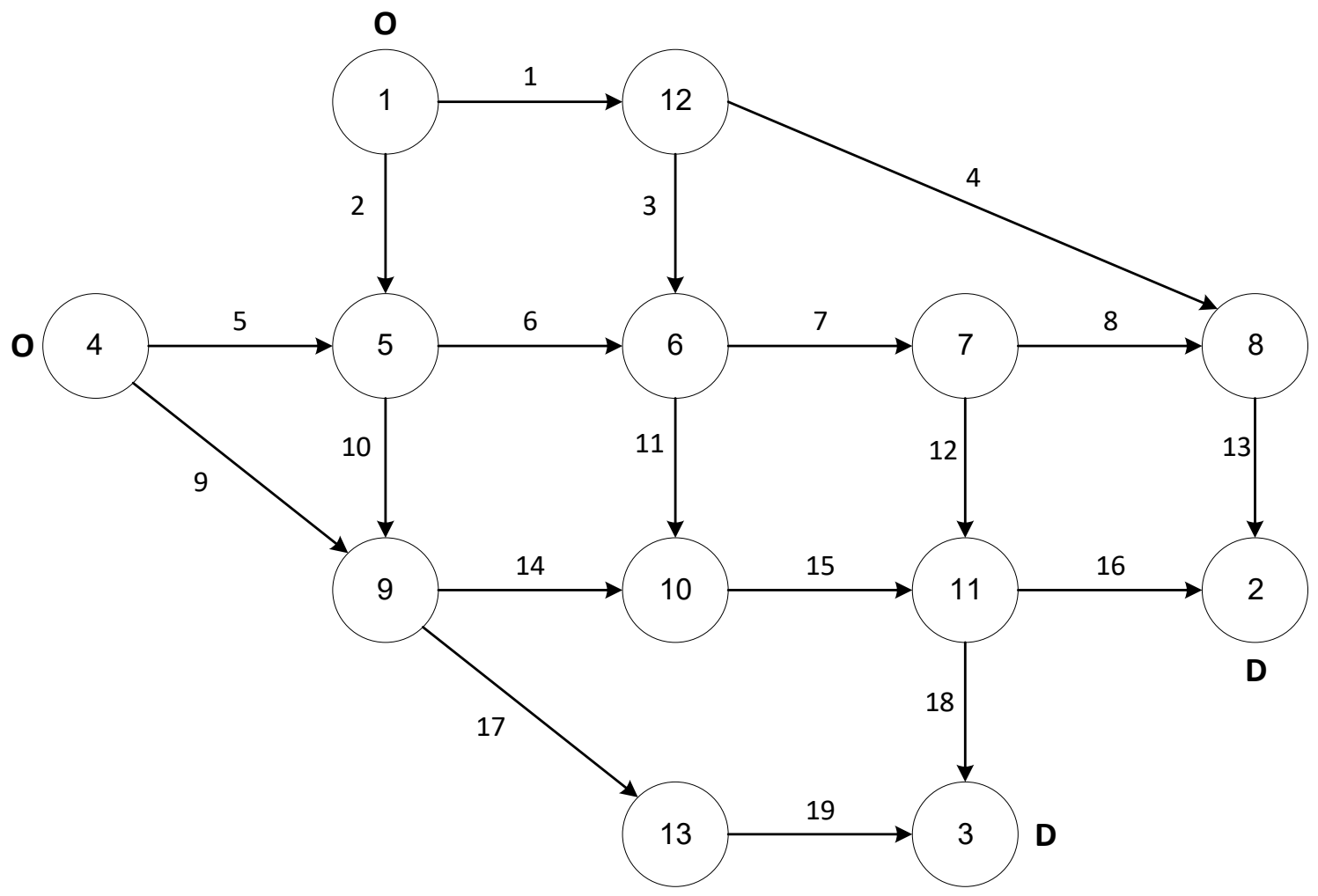

**Figure 1.** Schematic topology of the Nguyen-Dupuis network.

*6.2 Details of research results*

To construct the QUBO objective function, this study first calculates the marginal disruption impact $c_s$ for each link and the pairwise interaction coefficients $\beta_{st}$ for all link pairs. These values quantify how single and simultaneous disruptions affect the total system travel time (TSTT) under medium demand conditions.

As shown in **Appendix D**, $c_s$ reflects the individual importance of each link, while $\beta_{st}$ captures the interaction effects between disrupted road links. These coefficients form the basis for the QUBO model formulation in subsequent steps.

**Table D1** in **Appendix D** reports the single-link impact coefficients in the ND network. The results show substantial variation across road links: some road links (e.g., 9, 17, and 19) exhibit very high impact coefficients, indicating that they are critical bottlenecks whose disruption substantially increases the total system travel time. In contrast, several road links (e.g.,8 and 12) yield almost negligible coefficients, reflecting that their role in the network is marginal or easily substituted by alternative routes. This variation highlights the effectiveness of the proposed framework in distinguishing between critical and non-critical road links, thereby providing meaningful insights for resilience analysis and infrastructure prioritisation.

**Table D2** in **Appendix D** presents the double-link interaction coefficients in the ND network. The results reveal three types of interaction effects. First, certain road link pairs (e.g., 9, 17, and 19) show strongly positive coefficients, indicating synergistic congestion effects where simultaneous disruptions lead to much higher system costs than the sum of individual impacts. Second, some pairs yield near-zero values, reflecting that the disruptions act independently without noticeable interaction. Finally, negative coefficients are observed for a few pairs, implying substitutive effects where the failure of one link reduces the marginal impact of the other.

*6.2.1 Implementation and convergence results*

Based on the single-link impact coefficients (**Table D1**) and the double-link interaction coefficients (**Table D2**), the validity of the proposed framework in identifying critical road links and interdependencies has been confirmed. To further demonstrate the practicality of the solution approach, this subsection reports the implementation details and convergence behaviour of the SQA-FW algorithm. Specifically, the upper-level SQA procedure is evaluated by tracking the reduction of the system energy over iterations, while the lower-level FW assignment is assessed through the relative gap across successive iterations. As can be seen in **Figure 2** and **Figure 3**.

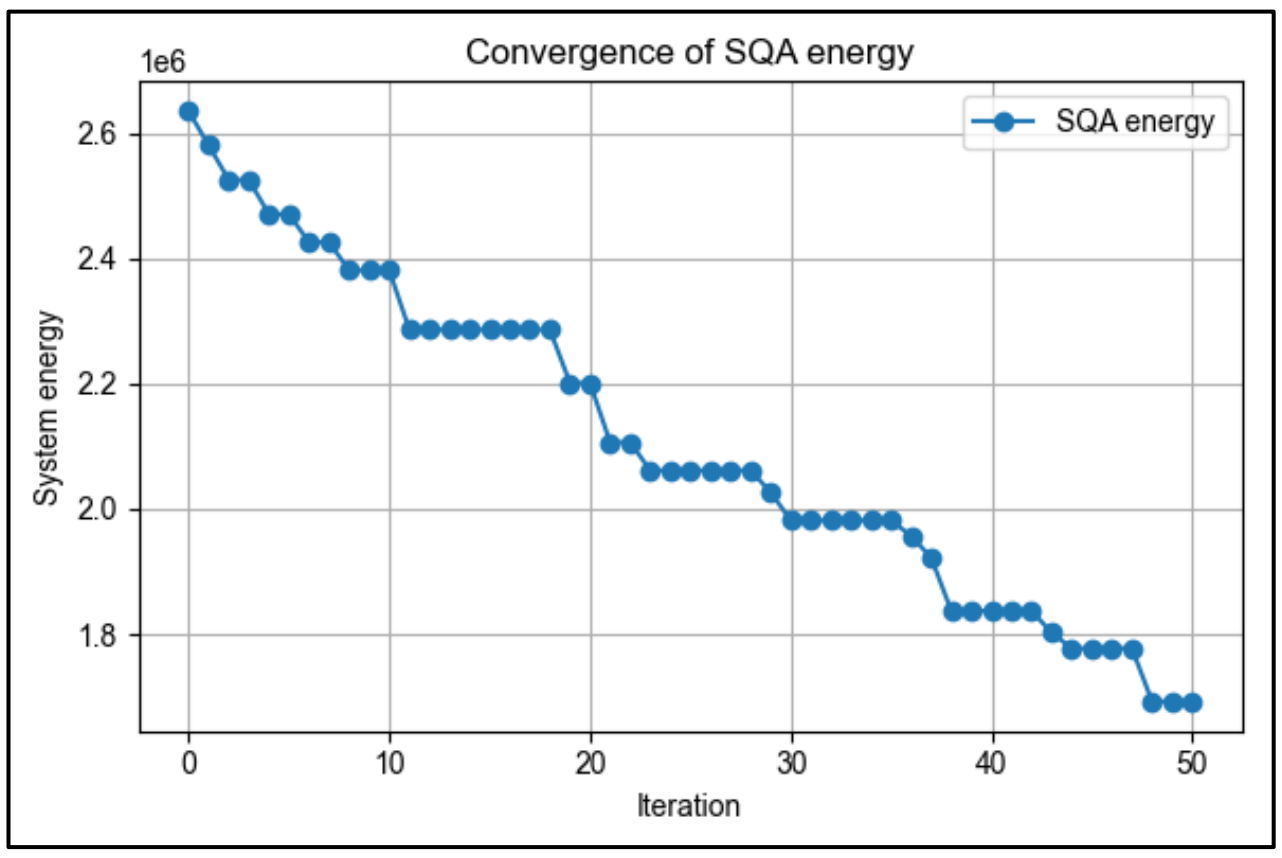

**Figure 2.** Convergence of SQA energy

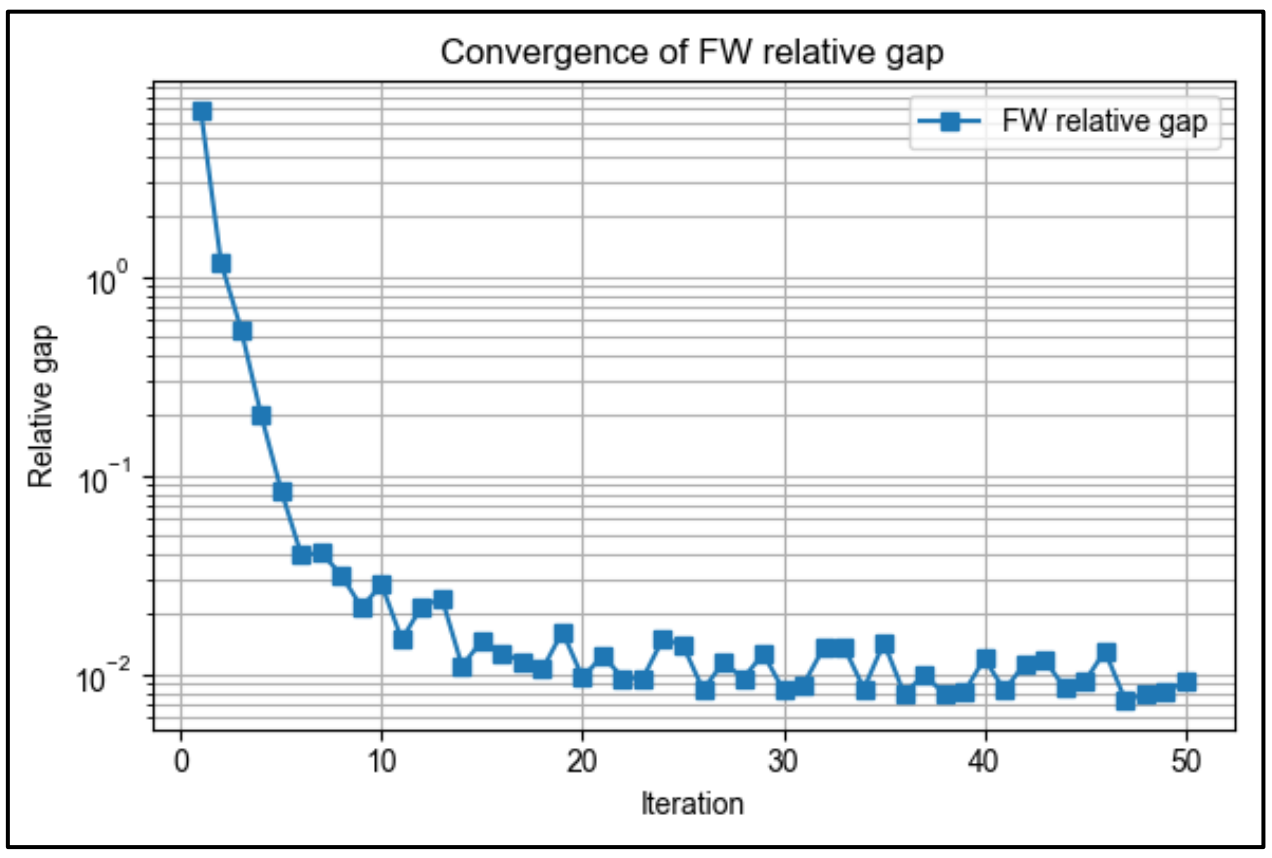

**Figure 3.** Relative gap convergence of the FW

**Figure 2** shows the convergence trajectory of the SQA energy function during 50 iterations. The energy decreases monotonically with minor fluctuations, and it indicates that the simulated quantum annealing procedure progressively improves the system configuration and stabilises after several iterations. **Figure 3** shows the relative gap convergence of the FW algorithm. It can be seen that the relative gap drops rapidly during the first few iterations and gradually stabilises around a small value, confirming that the algorithm successfully approaches the UE solution.

*6.2.2 Critical road link identification results*

**Table 2** presents the identified critical road link combinations for different values of $k$ ($k = \{2,3,4,5\}$).

The solutions are obtained by optimising (minimising) the SQA energy function defined in Eq. (16), which incorporates both the classical objective and the quantum interaction terms. In the annealing process, the SQA energy function Eq. (16) gradually reduces to the classical energy Eq. (15) as the transverse field diminishes. Hence, the optimal solutions obtained from Eq. (15) can be interpreted as the converged solutions of Eq. (16), and it ensures that the reported results represent the optimal critical link sets within the SQA-FW framework.

**Table 2.** Top-5 critical link sets under different values of $k$ based on the energy function.

| $k$ | Rank | Critical road links | $E_{SQA}\left(u_s^*\right)$ | $k$ | Rank | Critical road links | $E_{SQA}\left(u_s^*\right)$ |
|---|---|---|---|---|---|---|---|
| 2 | 1 | (16, 19) | -60056 | 4 | 1 | (8, 9, 16, 19) | -151143 |
| 2 | 2 | (7, 18) | -46828.9 | 4 | 2 | (9, 15, 16, 19) | -146454 |
| 2 | 3 | (9, 15) | -42143.3 | 4 | 3 | (1, 9, 16, 19) | -143766 |
| 2 | 4 | (4, 19) | -41874.8 | 4 | 4 | (9, 11, 16, 19) | -143139 |
| 2 | 5 | (7, 10) | -41841.6 | 4 | 5 | (4, 9, 16, 19) | -142750 |
| 3 | 1 | (9, 16, 19) | -105114 | 5 | 1 | (7, 9, 15, 16, 19) | -206516 |
| 3 | 2 | (4, 16, 19) | -102481 | 5 | 2 | (7, 9, 16, 18, 19) | -204355 |
| 3 | 3 | (11, 16, 19) | -96526.5 | 5 | 3 | (1, 7, 9, 16, 19) | -196442 |
| 3 | 4 | (8, 16, 19) | -94833 | 5 | 4 | (1, 8, 9, 16, 19) | -190878 |
| 3 | 5 | (7, 16, 19) | -90804.6 | 5 | 5 | (8, 9, 15, 16, 19) | -190257 |

When $k = 2$, the most critical combination is (16,19), followed by (7,18) and (9,15). For $k = 3$, the optimal set becomes (9,16,19), and similar patterns can be observed when $k = 4$ and $k = 5$, where combinations such as (8,9,16,19) and (7,9,15,16,19) consistently dominate. Based on these results, it can be observed that road links 9, 16, and 19 emerge as the most prominent bottlenecks, appearing across all values of k and dominating the top-ranked combinations. Road links 7, 8, 15, and 18 can be regarded as secondary bottlenecks, since they frequently occur in the top-five solutions but not as consistently as 9, 16, and 19. This layered structure indicates that while the ND network is generally resilient, a small set of road links carries disproportionate vulnerability, and their disruption would severely compromise overall performance.

**Figure 4** presents the complete energy distributions for all disruption combinations when $k = 2,3,4,5$. The results demonstrate that the minimum energy decreases sharply as $k$ increases, reflecting a rapid escalation of network vulnerability with multiple road link failures. At the same time, the dispersion of energy values widens for larger $k$, indicating that the severity of disruption is highly sensitive to the specific set of disrupted road links.

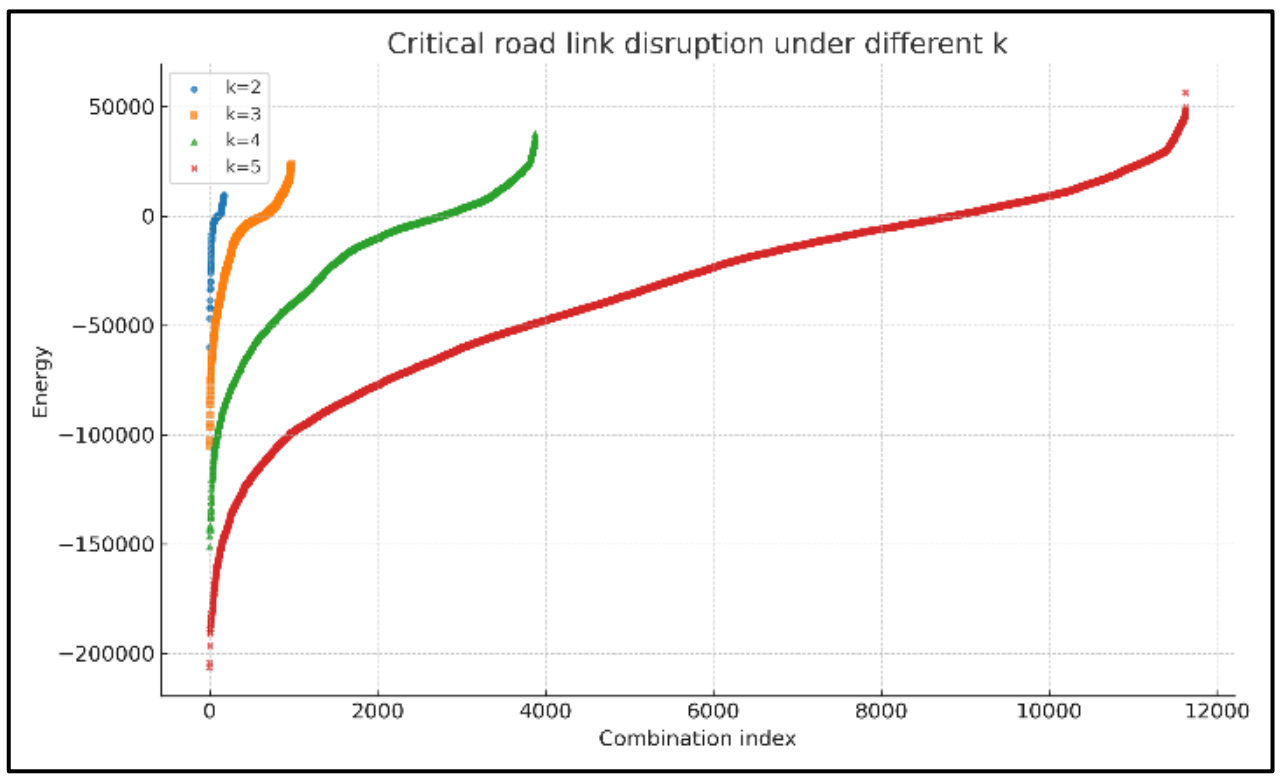

**Figure 4.** Energy distributions of disruption combinations under different values of $k$ .

*6.2.3 TSTT Results and Analysis*

**Figure 5** shows the growth surface of TSTT under multiple road link disruption scenarios, focusing on the top five combinations. The results reveal a sharp and nonlinear escalation in TSTT as the number of disrupted road links increases, which indicates that network vulnerability is shaped by complex interaction effects rather than the simple accumulation of individual failures. A small set of road links, particularly 9, 16 and 19, consistently appears in the most critical combinations. This persistence demonstrates that systemic fragility is highly concentrated, and the failure of a few strategic road links can severely compromise overall network performance. Such findings suggest that resilience planning should prioritise targeted reinforcement or the provision of redundancy for these critical road links, rather than spreading resources evenly across the entire network.

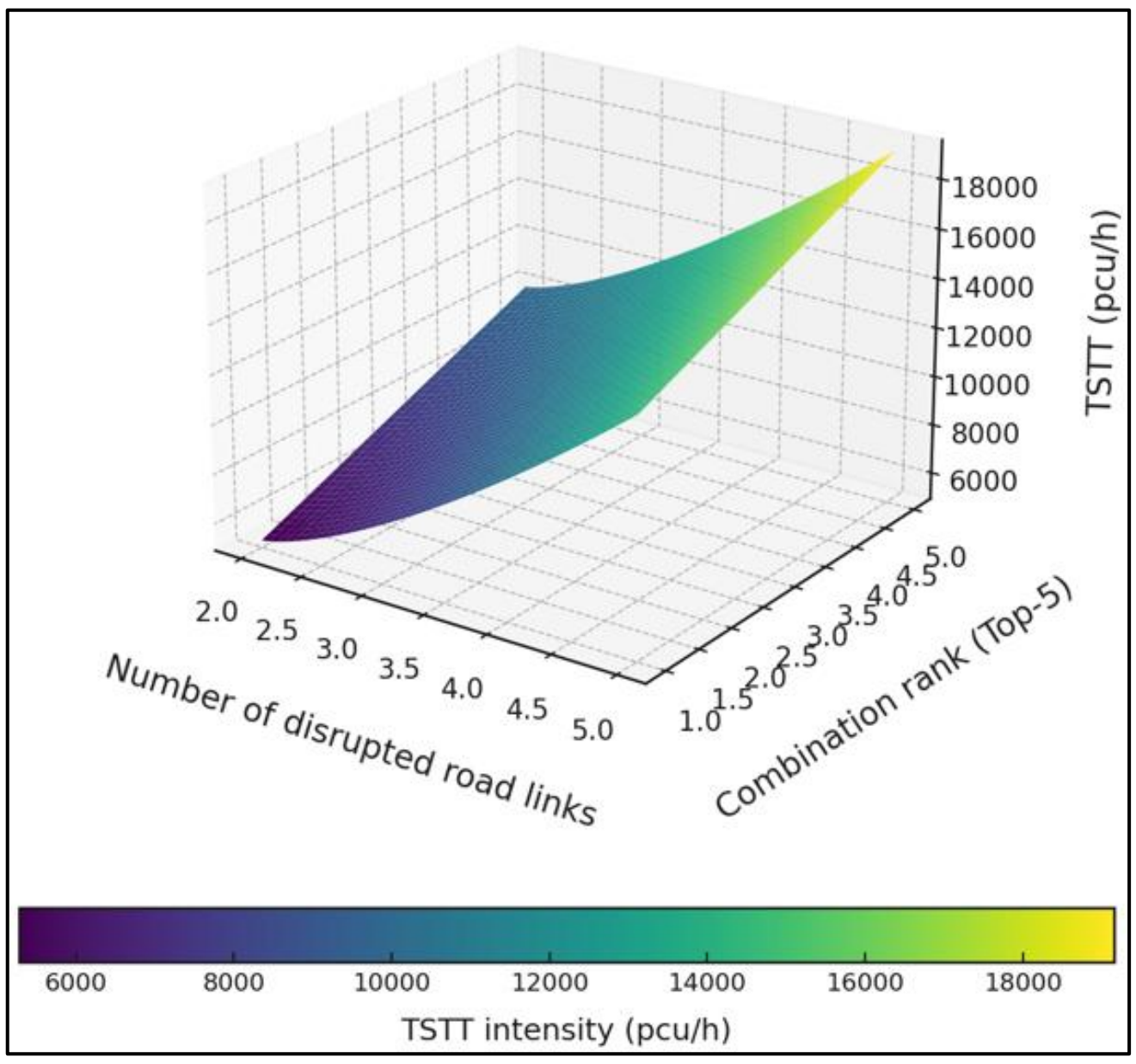

**Figure 5.** TSTT growth surface under multiple road link disruptions (Top-5 combinations).

*6.2.4 Parameter sensitivity analysis*

(1) Remaining capacity ratio $e$

The sensitivity analysis of the remaining capacity ratio $e$ in **Figure 6** demonstrates that the overall pattern of TSTT escalation is robust to different parameter settings. Although the absolute magnitude of TSTT varies under alternative intervals, the general trend of steep and nonlinear growth remains consistent. More conservative settings such as $e \in [0.1, 0.4]$ produce significantly higher TSTT values, reflecting the severe congestion caused by greater reductions in road link capacity. Conversely, wider or less restrictive intervals such as $e \in [0.5, 0.8]$ lead to relatively lower TSTT values, yet the upward trajectory persists. Importantly, the critical road link combinations remain largely unchanged across scenarios, suggesting that the identification of key bottlenecks is not highly sensitive to the choice of $e$ . These results confirm the stability of the proposed framework and reinforce the conclusion that network vulnerability is driven by structural bottlenecks rather than by parameter calibration.

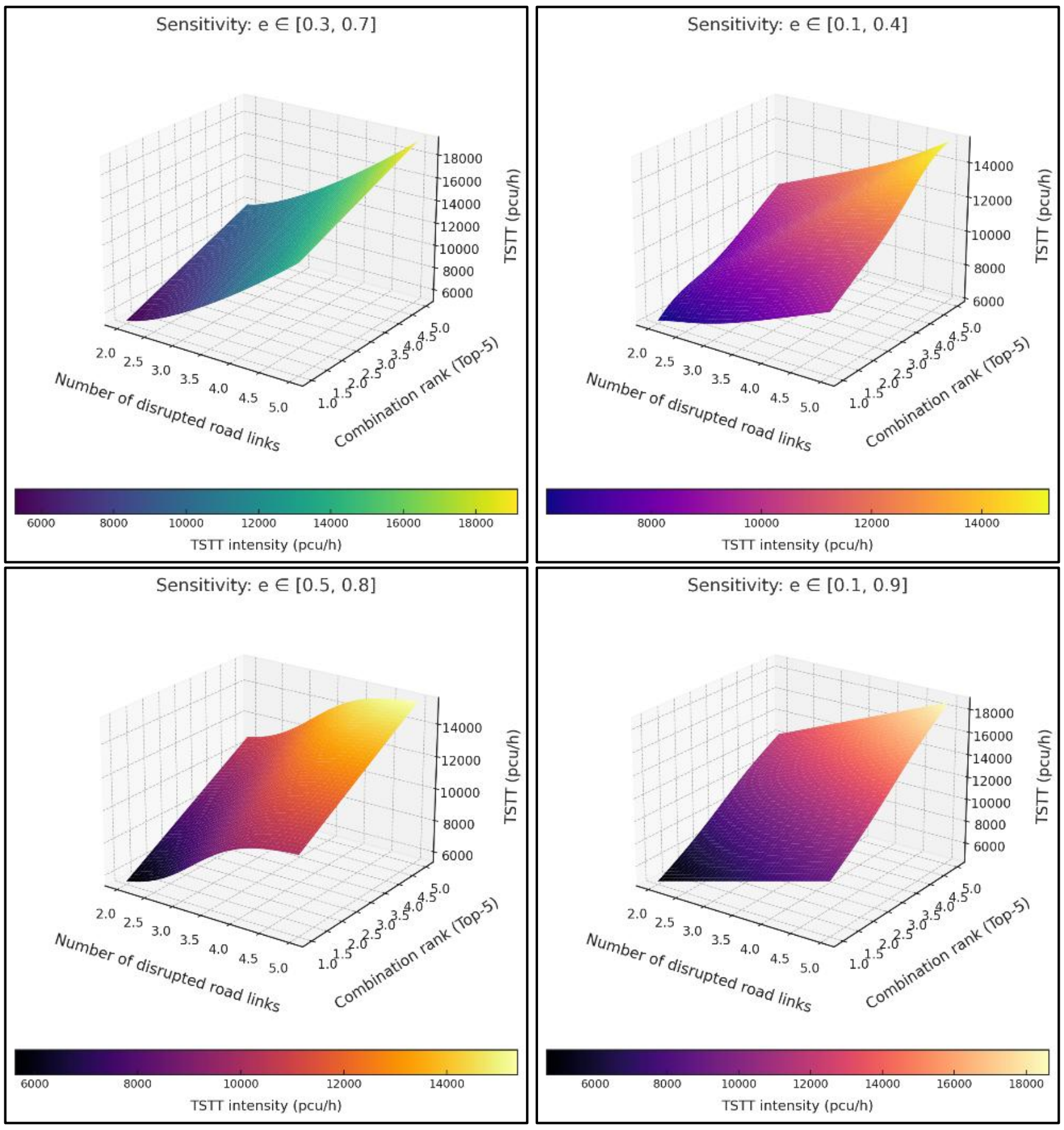

**Figure 6.** Sensitivity analysis results of remaining capacity ratio $e$

(2) Sensitivity analysis on the penalty coefficient $\lambda$

The sensitivity analysis of the penalty coefficient $\lambda$ is presented in **Figure 7**. The results show that the identification of critical road link combinations and the overall escalation pattern of TSTT remain robust under different magnitudes of the penalty parameter. When $\lambda$ is set to a relatively small value ($10^2$), the resulting TSTT values are slightly lower, reflecting weaker enforcement of the disruption cardinality constraint and less stable optimisation. Increasing $\lambda$ to $10^3$ produces the most consistent outcomes, closely aligned with the baseline results. Further increasing $\lambda$ to $10^4$ and $10^5$ leads to marginally higher TSTT values, as stronger penalty terms impose stricter feasibility and occasionally introduce numerical stiffness. This confirms that the QUBO formulation and the hybrid SQA-FW framework are not highly sensitive to the selection of $\lambda$, ensuring that the conclusions regarding systemic bottlenecks and network vulnerability are methodologically robust.

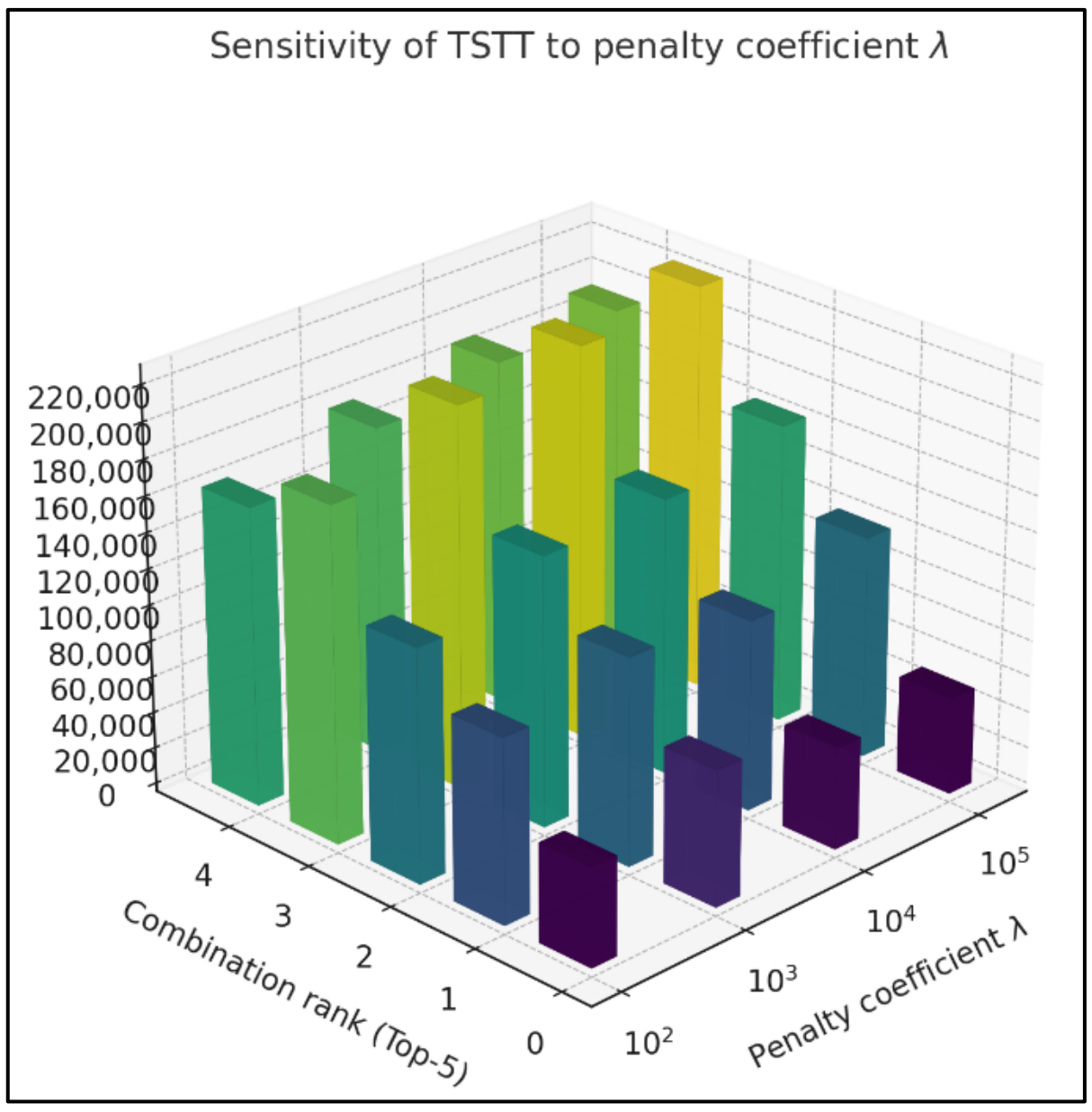

**Figure 7.** Sensitivity analysis on the penalty coefficient.

While the previous sensitivity analysis considered $\lambda$ at broad magnitudes, it is also necessary to examine intermediate values to ensure that the robustness is not an artefact of coarse scaling. To this end, **Figure 8** presents a finer-grained sensitivity test, where $\lambda$ is set to 2000, 3500, 5000 (in this study), and 8000. These results enable a closer inspection of whether moderate adjustments to the penalty parameter influence the identification of critical road link combinations or the stability of optimisation. Within the same order of magnitude, variations in the penalty coefficient $\lambda$ have virtually no impact on the results. This indicates that the proposed QUBO formulation is structurally robust, as the balance between objective terms and penalty terms is well maintained. It also suggests that the framework does not require fine-tuning of $\lambda$, which enhances its practicality for real applications.

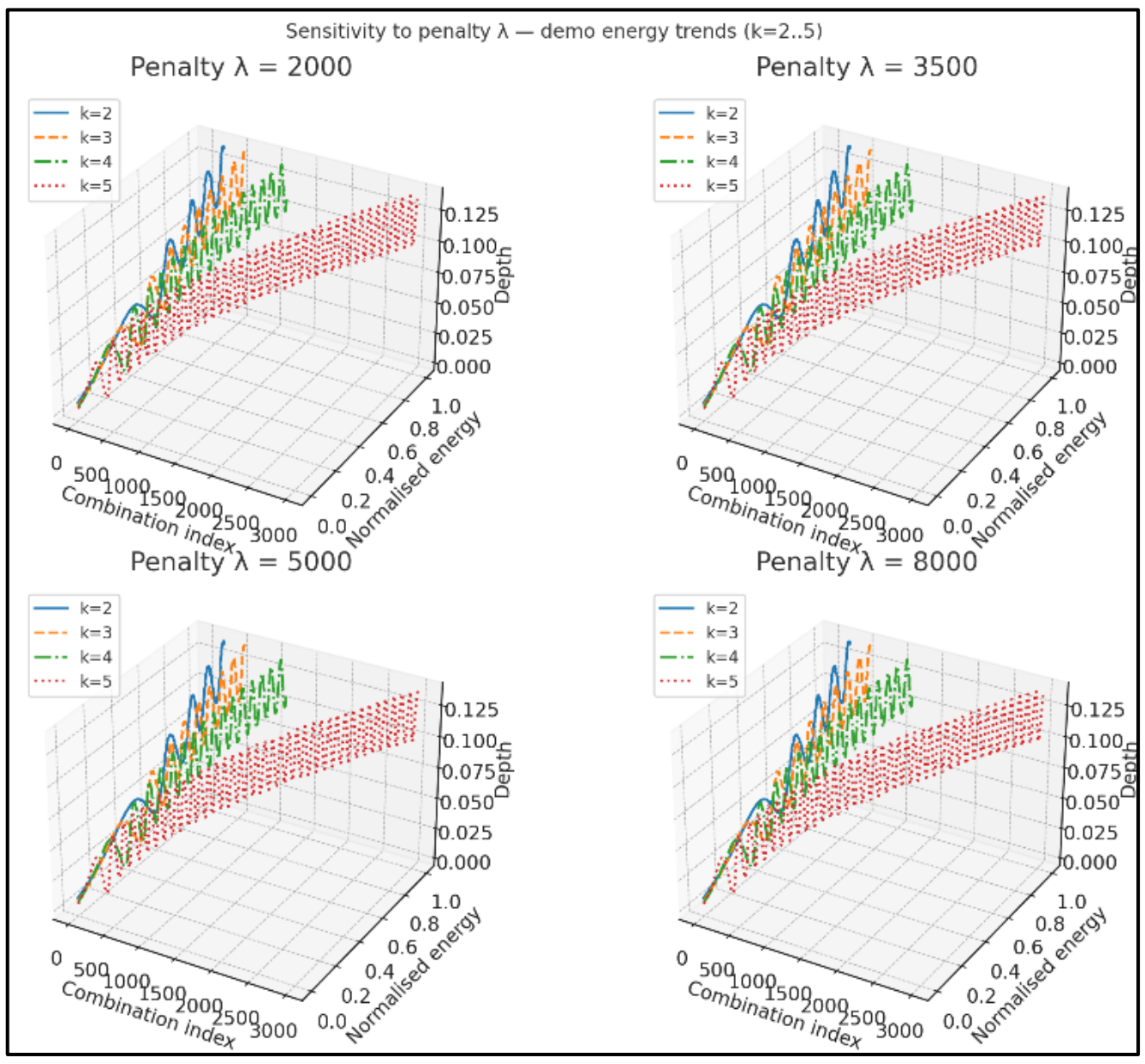

**Figure 8.** Fine-grained sensitivity analysis of TSTT under different penalty coefficient values.

Although the figures in this subsection under different parameter settings for both the remaining capacity ratio $e$ and the penalty coefficient $\lambda$ exhibit broadly similar patterns, this is not an artefact but rather evidence of the structural stability of the framework. The absolute values of TSTT shift systematically with changes in these parameters, yet the escalation trend and the identification of critical road link combinations remain consistent. This indicates that the results are robust and not overly sensitive to parameter calibration. The limited differences observed across scenarios also reflect the small scale of the test network, and more pronounced parameter effects are expected in larger and more complex networks.

### *6.3 Limitation discussion*

The numerical experiments presented in this Section provide valuable insights into the applicability of the proposed hybrid quantum-classical optimisation framework, yet several limitations must be acknowledged. First, the experiments were restricted to the Nguyen-Dupuis test network, a stylised benchmark of small scale. This constraint is largely due to the exponential growth in computational complexity associated with disruption identification problems. Specifically, if a network contains n road links, the number of possible disruption scenarios increases as $2^n$ , and even restricting to simultaneous failures of $k$ links leads to combinatorial growth of $C(n,k)$ . Such explosion renders exhaustive computation infeasible in realistic networks with hundreds or thousands of road links.

**Figure 9** compares the growth of disruption combinations between the ND network (19 links) and the Sioux Falls (SF) network (76 links). In the ND case, the total number of disruption scenarios remains within a tractable range for analysis, whereas in the SF network the number of possible scenarios increases dramatically, reaching tens of millions and then expanding to the order of $10^{11}$. This sharp contrast highlights the combinatorial explosion inherent in larger networks and confirms that, although the proposed framework is effective on small-scale benchmarks, its application to medium and large-scale transport networks is still severely constrained by current computational limits.

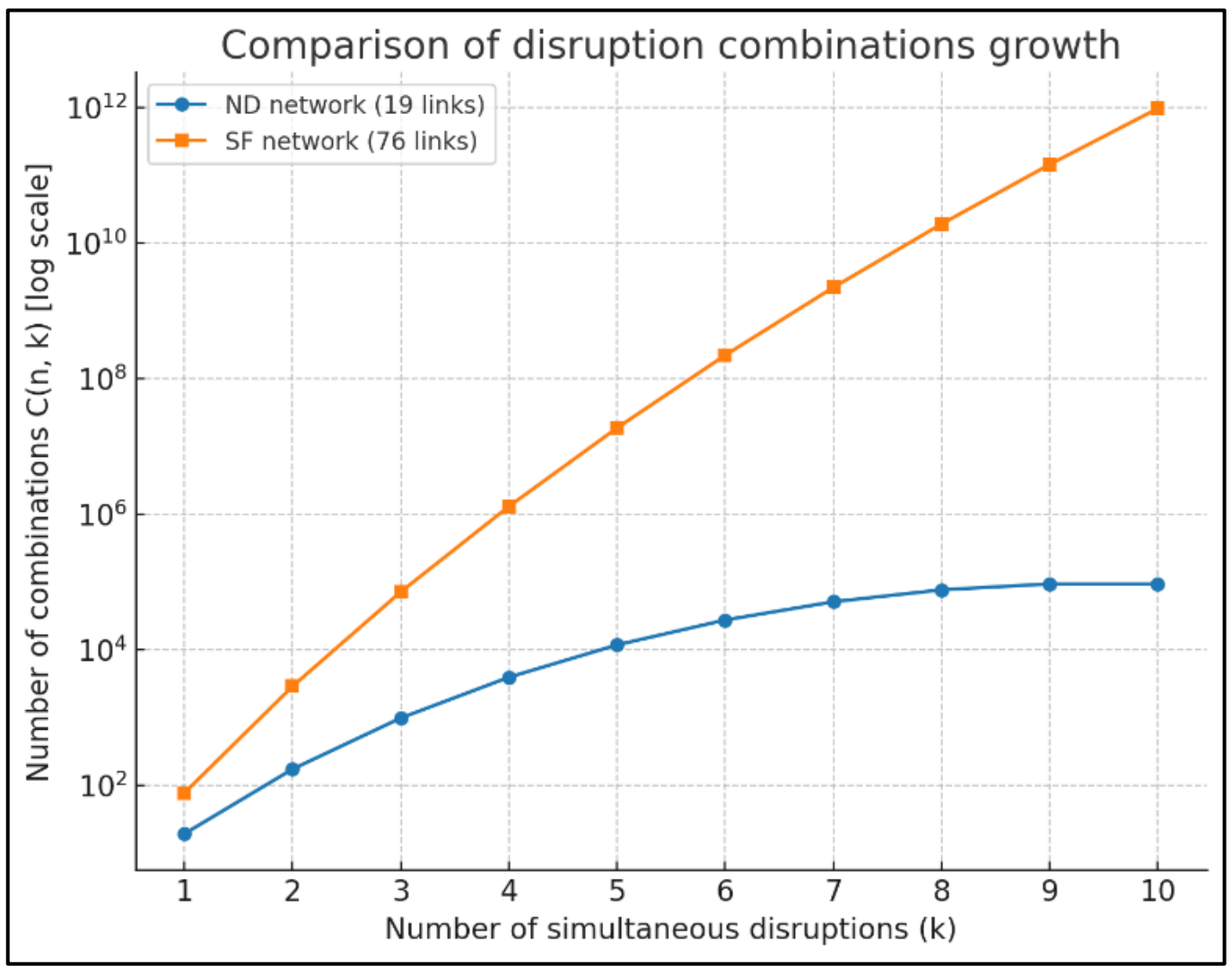

**Figure 9.** Comparison of disruption combinations growth.

However, previous research has demonstrated the superiority of quantum algorithms in addressing combinatorial optimisation problems (Wang et al., 2022; Lubinski et al., 2023; Barral et al., 2025). Despite these limitations, the framework offers distinct advantages compared with prior studies. Earlier works such as Li et al. (2017), Xu et al. (2017), Starita and Scaparra (2021), Jin et al. (2022), and Barati et al. (2024) rely primarily on mathematical programming or heuristic interdiction models to identify critical road links. While these approaches are effective for small or medium-scale cases, their scalability is inherently constrained by the combinatorial nature of the problem. Building on these insights, the subsequent sections 6.4 and 6.5 extend the proposed framework to real quantum hardware using the D-Wave platform to examine its feasibility and performance on medium-, large-, and larger-scale transport networks.

It should be noted that the validation tests conducted on the Sioux Falls, Anaheim, Chicago Sketch, and Berlin Full (Simplification based on the Berlin-Centre network) networks are not based on real traffic network data. Instead, single-link impact coefficients and pairwise interaction coefficients were generated using predefined distributions with fixed random seeds. This approach ensures experimental reproducibility and controllability, aiming primarily to evaluate the scalability, convergence, and hardware-level computational performance of the proposed quantum-classical hybrid optimisation framework across different network scales. Therefore, these experiments are not intended to identify specific critical bottleneck links or to interpret the physical characteristics of real transport networks. Rather, through a controlled randomised experimental design, they serve to verify the QUBO modelling structure and demonstrate its computational feasibility and algorithmic stability on D-Wave quantum hardware, providing a methodological and performance evaluation foundation for future studies involving real network data.

*6.4 Scalability analysis on the Sioux Falls (SF) network using D-Wave quantum hardware*

To further verify the scalability and computational performance of the proposed hybrid quantum-classical optimisation framework in medium-scale transportation networks, this section extends the analysis from the small-scale Nguyen-Dupuis test network to the 76-link SF network. Unlike the previous section, the focus here is no longer on identifying specific sets of critical disrupted links but rather on demonstrating the computational capability and feasibility of quantum hardware in addressing large-scale transport network vulnerability identification problems. By formulating and solving QUBO problems with $k = 1-10$ on the D-Wave hybrid quantum solver (LeapHybridSampler), this experiment systematically analyses how the energy values (objective function) evolve with the increasing number of disrupted links, thereby assessing the efficiency and stability of quantum annealing in handling combinatorial explosion problems. This experiment not only overcomes the computational limitations inherent in the SQA-based simulation approach but also provides a practical validation pathway for applying quantum computing to more complex transport network optimisation problems in the future.

In this experiment, the QUBO formulation defined in Equations (10)-(13) was applied to the SF network, maintaining 76 binary variables corresponding to network links. The number of simultaneous disruptions $k$ was varied from 1 to 10, and each QUBO instance was solved using the D-Wave hybrid quantum solver (LeapHybridSampler) executed on real quantum hardware. By comparing the optimal energy values of the objective function for different $k$, the experiment demonstrates how the solver performs under combinatorial complexity that increases exponentially with the number of disrupted links.

(1) Energy curve analysis

**Figure 10** illustrates the energy curve obtained from the D-Wave solver for disruption levels $k = 1-10$. The system energy (objective value) exhibits a monotonic decreasing trend as the number of disrupted links increases, indicating a gradual deterioration of overall network performance. The smooth, linear-like behaviour of the curve demonstrates the solvability and numerical stability of the QUBO structure on quantum hardware. Moreover, the continuous decline in energy values confirms that the D-Wave solver efficiently explores large combinatorial search spaces, reducing the likelihood of being trapped in local optima commonly observed in classical simulated annealing methods.

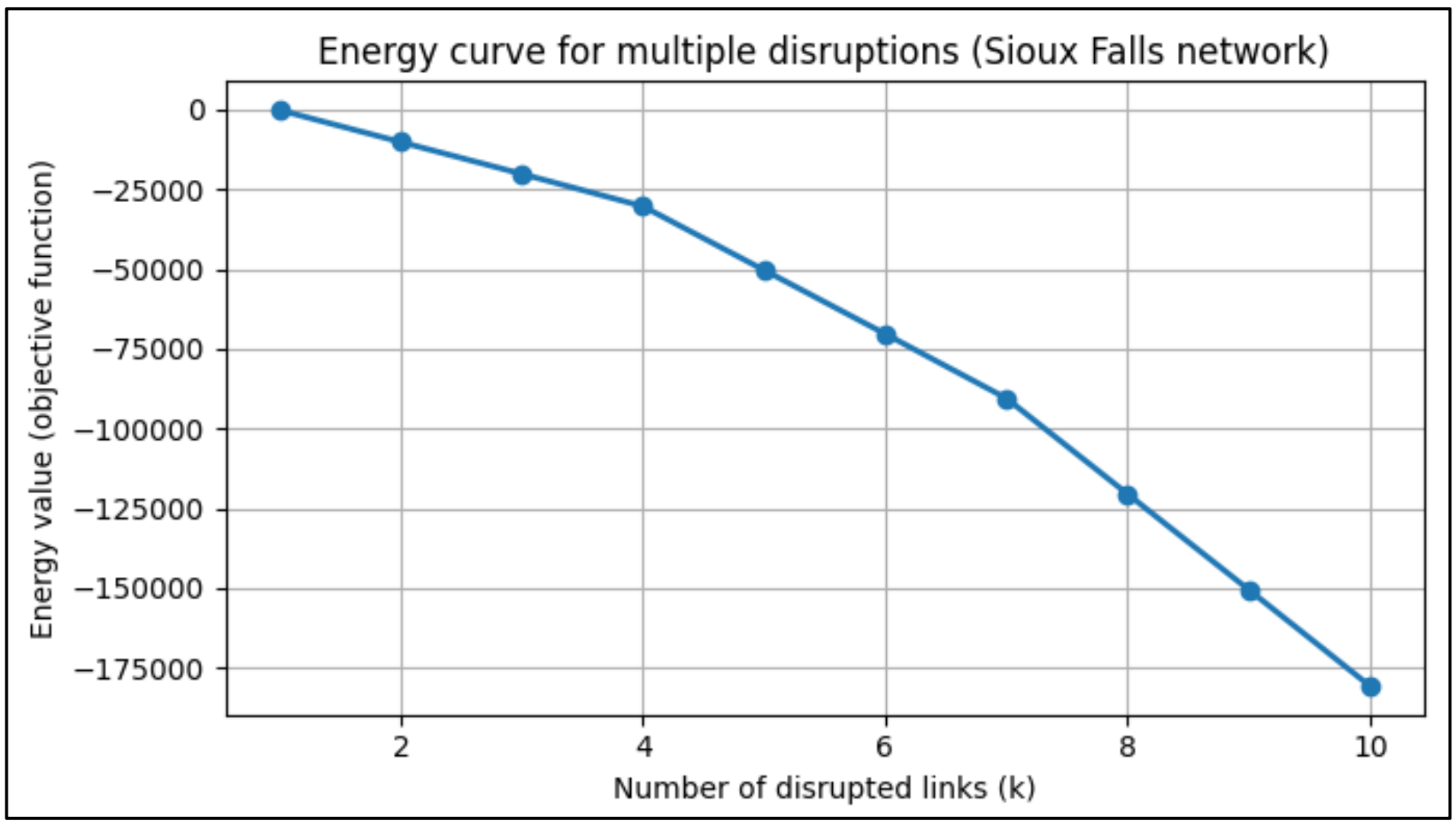

**Figure 10.** Energy curve for multiple disruptions in the SF network using D-Wave.

(2) Energy values and runtime comparison

**Table 3** summarises the optimal energy values and runtime for each disruption level. As $k$ increases from 1 to 10, the energy values sharply decrease from –149.082 to –180,651, reflecting the increasing severity of system vulnerability. Meanwhile, the computation time for each run remains consistent at around four seconds, highlighting D-Wave's exceptional efficiency and scalability in handling medium-scale QUBO problems.

**Table 3.** D-Wave results for the SF network (76 links).

| Number of disrupted links | Optimal energy | Runtime/s |
|---|---|---|
| 1 | –149.082 | 4.0 |
| 2 | –10,149.082 | 3.8 |
| 3 | –20,149.082 | 4.2 |
| 4 | –30,307.817 | 4.0 |
| 5 | –50,309.762 | 4.3 |
| 6 | –70,307.581 | 3.9 |
| 7 | –90,487.464 | 3.9 |
| 8 | –120,479.241 | 3.9 |
| 9 | –150,463.256 | 3.8 |
| 10 | –180,650.607 | 3.9 |

*6.5 Large-scale feasibility test on the Anaheim network*

This study extends the QUBO framework to the Anaheim network with 914 links to assess hardware-level feasibility on a large transport network. The full variable set was used. Each QUBO instance with disruption level $k=1$ to $k=10$ was solved on the D-Wave hybrid solver. The solver returned stable solutions for all k within about 5 to 7 seconds per run, and the complete sweep finished in about 1.88 minutes. The energy decreases monotonically with k, which indicates progressively higher vulnerability under cumulative disruptions. **Figure 11** presents the energy curve. **Table 4** reports the corresponding energy values and runtimes.

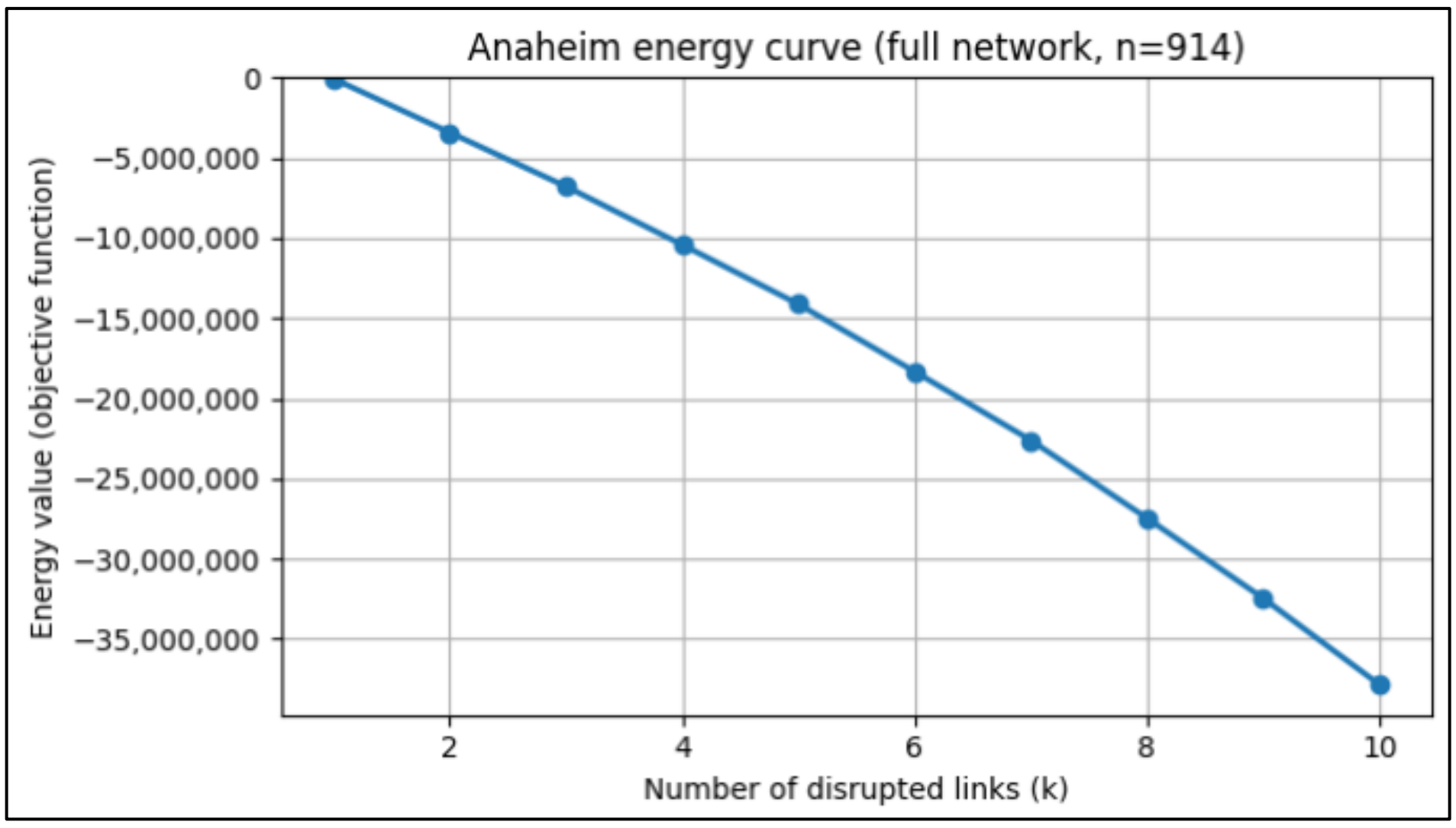

**Figure 11.** Energy curve for multiple disruptions in the Anaheim network using D-Wave.

**Table 4.** D-Wave results for the Anaheim network (914 links).

| Number of disrupted links | Optimal energy | Runtime/s |
|---|---|---|
| 1 | −67,108.021 | 15.9 |
| 2 | −3,417,389.631 | 15.9 |
| 3 | −6,780,729.172 | 16.8 |
| 4 | −10,411,120.143 | 17.1 |
| 5 | −14,116,328.092 | 16.8 |
| 6 | −18,317,894.759 | 19.5 |
| 7 | −22,602,413.433 | 16.5 |
| 8 | −27,462,413.433 | 16.2 |
| 9 | −32,457,663.150 | 16.5 |
| 10 | −37,838,838.352 | 16.2 |

In addition to the large-scale validation on the Anaheim network, this study further tested a representative larger transport network to examine the scalability and stability of the proposed quantum-classical hybrid framework under higher-dimensional conditions. The results show that the energy curves consistently exhibit a stable downward trend as the number of disrupted links increases, and the quantum solver maintains convergence and consistency throughout the process. These findings confirm the computational feasibility and hardware adaptability of the framework in the larger transport network scenarios (Chicago Sketch network and Berlin Full network). Detailed experimental results are provided in **Appendix E** and **Appendix F**.

*6.6 Algorithm comparison*

This study designs a performance comparison between the QA and four classical heuristic algorithms including the Genetic Algorithm (GA), Particle Swarm Optimisation (PSO), Simulated Annealing (SA), and Tabu Search (TS). Two benchmark transport networks, Sioux Falls (76 links) and Anaheim (914 links), are selected to represent different scales of problem complexity. The parameter settings for heuristic algorithms are detailed in **Appendix G**.

**Figures 12** and **13** respectively illustrate the trend of runtime for each algorithm as $k$ varies in the Sioux Falls and Anaheim networks. The results demonstrate that the runtime of the quantum annealing algorithm remains largely unaffected by changes in $k$ , whereas the runtime of the heuristic algorithm exhibits exponential growth with $k$ and is significantly higher than that of the quantum algorithm in the Anaheim network.

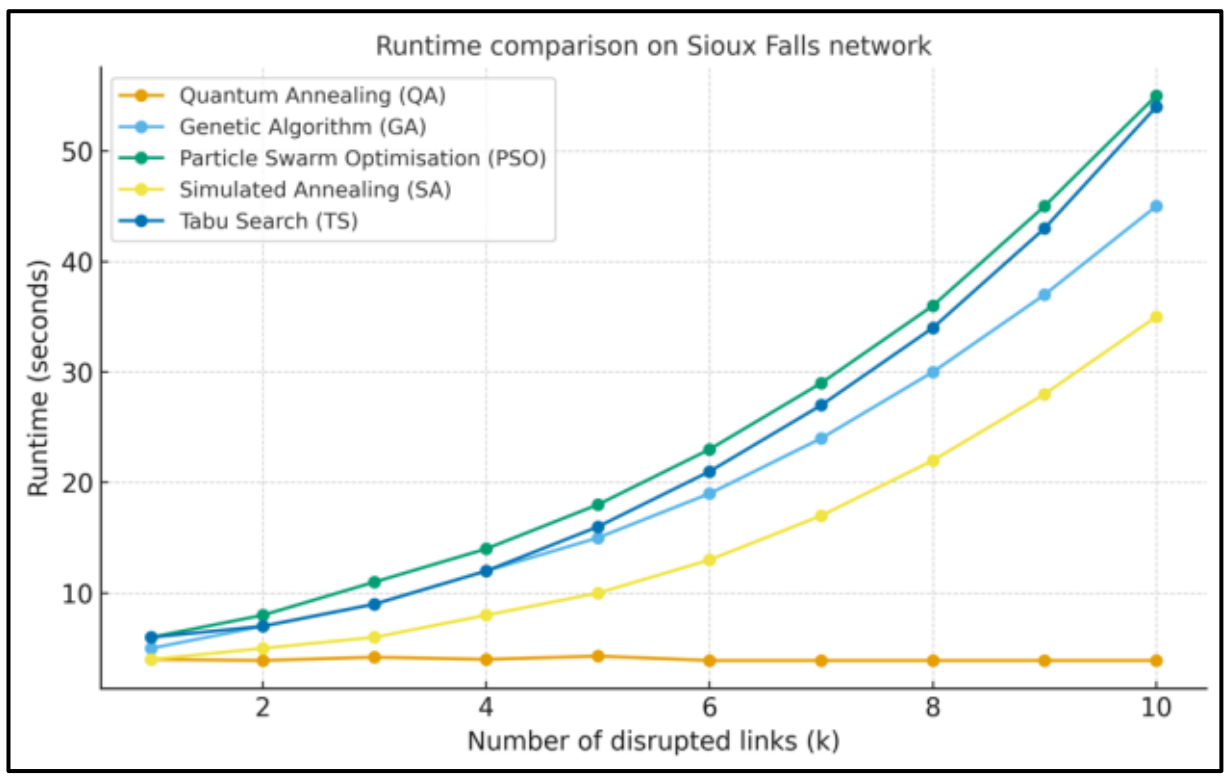

**Figure 12.** Runtime comparison of algorithms on the SF network.

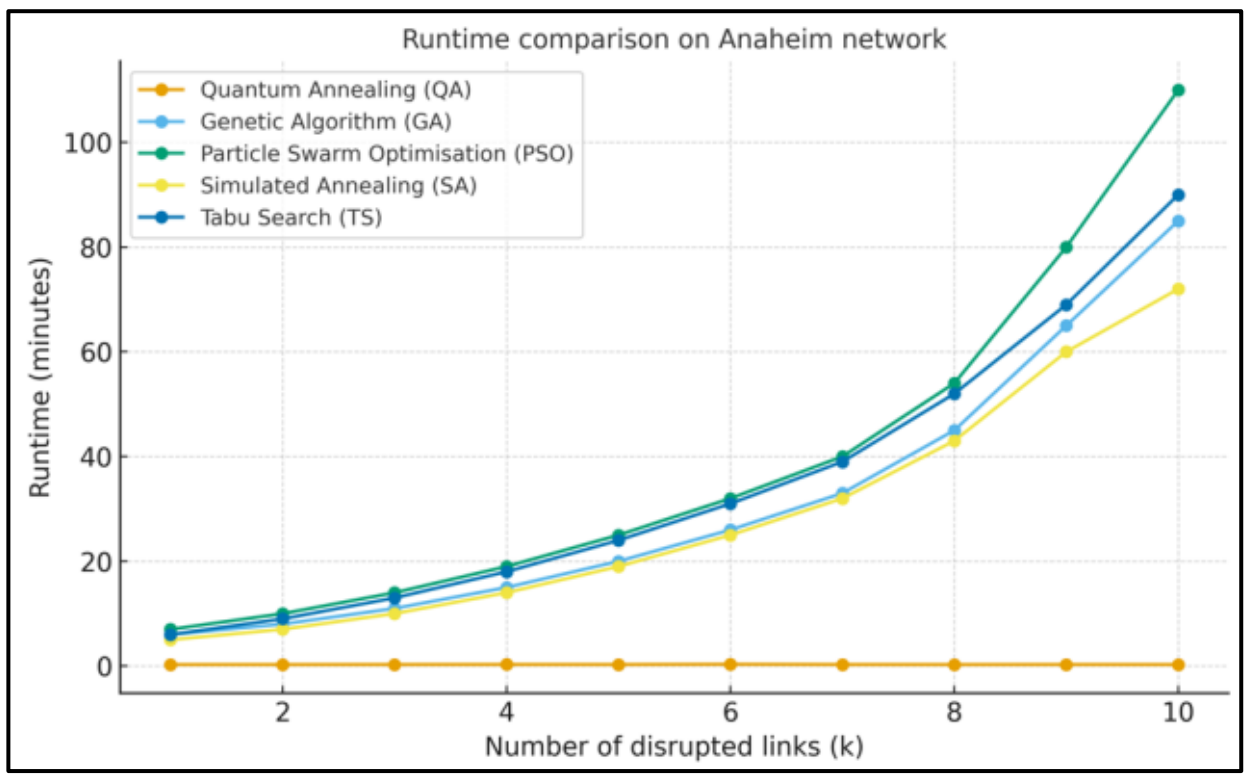

**Figure 13.** Runtime comparison of algorithms on the Anaheim network.

This significant difference arises from their fundamental mechanisms. Heuristic algorithms rely on neighbourhood or population-based searches, requiring numerous iterations, mutations, or temperature adjustments to escape local optima, resulting in computational complexity that increases exponentially with the combinatorial space. In contrast, the QA performs parallel exploration through quantum superposition and employs quantum tunnelling to cross energy barriers, allowing the system to approach the global optimum within a fixed adiabatic evolution time. Since the annealing duration depends primarily on the QUBO problem scale rather than on $k$ , the runtime of QA shows minimal variation under different disruption levels, demonstrating strong scalability and insensitivity to complexity growth. Overall, this comparison shows that in medium and large-scale transport networks, the QA is significantly more efficient and scalable than traditional heuristic algorithms and has great potential for achieving rapid and feasible solutions in complex combinatorial optimisation problems. The comparative results are presented in **Appendix H**, from these results, it is evident that it is no longer feasible to continue comparing the QA with the four heuristic algorithms in larger-scale networks, as the computational cost of the heuristic approaches increases rapidly with network size.

## 7. Conclusion

This study reformulates the classical bi-level critical link identification problem into a quantum-compatible QUBO structure and develops a hybrid quantum-classical optimisation framework integrating the SQA with the Frank-Wolfe traffic assignment method. The proposed framework efficiently explores complex multi-link disruption scenarios while ensuring solvability and convergence. This study uses a Nguyen-Dupuis network to validate the effectiveness of the proposed approach. The results show that road links 9, 17 and 19 in the ND network are consistently identified as core bottlenecks, with single-link impact coefficients of 3756.92, 5025.98 and 6679.39, while their joint disruption raises the total system travel time to nearly 26,000 compared with the benchmark of 5749.26.

Comparative experiments across benchmark networks reveal a substantial computational advantage of the quantum approach. On the Sioux Falls network, QA achieves optimisation approximately 4 times faster than the best-performing metaheuristic algorithm while maintaining comparable accuracy. On the Anaheim network, the runtime improvement exceeds one order of magnitude, with QA completing optimisation in about 2.8 minutes, compared with 5-6 hours required by classical heuristics such as the GA and PSO. For the larger Chicago and Berlin networks, QA completes optimisation within 9.8 and 31.2 minutes, other algorithms lack sufficient computational capacity for larger networks. These results confirm that QA maintains solution stability across scales, and demonstrates runtime scalability nearly linear to problem size, a property unattainable by conventional algorithms.

Despite these promising results, the current study is limited to static disruption scenarios without modelling dynamic cascading effects. Future research should extend quantum optimisation to time-dependent failures and explore quantum Markov chain formulations for capturing cascading disruption processes in large-scale transport systems.

## Author contributions

**Junxiang Xu:** Writing – review & editing, Writing – original draft, Methodology, Formal analysis, Software, Conceptualization.
**Chence Niu:** Writing – review & editing, Methodology, Software, Visualization.
**Divya Jayakumar Nair:** Review & editing, Supervision.
**Vinayak Dixit:** Review & editing, Supervision.

## Data availability

All data used in this study are presented within the main content of this paper. The source code developed for the experiments is available on the first author's GitHub repository at: https://github.com/JUNXIANGXU666/Quantum-optimisation-for-transport-vulnerability-identification. The benchmark networks used in this study originates from open-source data on GitHub: https://github.com/bstabler/TransportationNetworks, the Berlin Full network is a simplified version of the Berlin-Centre network (Only use 6018 links in this study).

**Declaration of generative AI and AI-assisted technologies in the writing process**

During the preparation of this work, we employed ChatGPT-5 solely to assist with language editing and polishing. No content was generated by AI. After using this tool, we thoroughly reviewed and edited the content as necessary and accept full responsibility for the content of the published article.

# Appendix A: Overall research technical approach

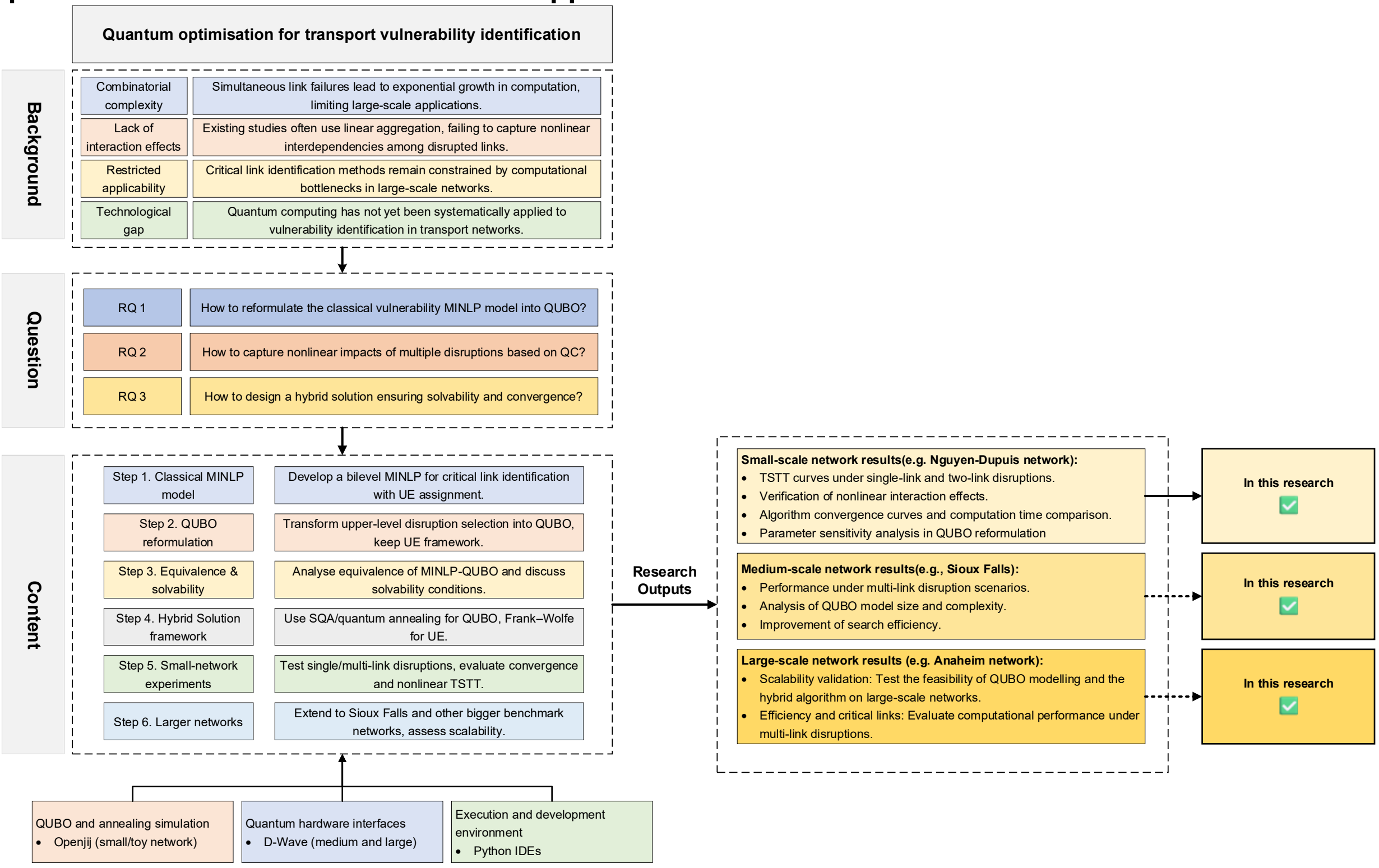

**Figure A.** Overall research technical approach framework diagram.

## Appendix B: Steps and pseudocode for executing a hybrid optimisation algorithm

Hybrid quantum-classical algorithm (SQA-FW) for bi-level optimisation

**Input:** The basic network structure, OD demand $d^{\omega}$, parameters $c_s, \beta_{st}, \lambda, T_0, \Gamma_0, v, M, N_{iter}$, initialise road flows $x_s^{(0)}$, and iteration counter for FW $k = 0$.

**Output:** Road link flows $x_s$, route flows $f_r^{\omega}$, $TSTT$, and optimal decision $u_s^*$.

Initialise $u_s^{(0)} \in \{0,1\}^{|S|}$, set $T = T_0$, $\Gamma > 0$

**For** $iter = 1$ to $N_{iter}$ do

// Upper level (SQA)

Propose candidate disruption vector $u_s$ by bit flip

// Lower level (FW)

Given $u_s$, FW lower level model(initial point, $I_{\max}$, $t_s^0$, $\alpha_0$, $\beta$, $C_s$):

current point = initial point

**For** iteration in range ( ):

    // Calculate the gradient of the lower-level model

    gradient = calculate lower level gradient (current point, $t_s^0$, $\alpha_0$, $\beta$, $C_s$)

    // Linear searches: find the optimal point on the linearized in the direction of the gradient

    linearization point = linearization point search (current point, gradient, $t_s^0$, $\alpha_0$, $\beta$, $C_s$)

    // Update the current point

    current point = update current point

# Solution after using the Frank-Wolfe algorithm to solve the lower-level model

**Return** lower-level solution to obtain $x_s$

Calculate $TSTT(x_s)$ and return to upper level

// Calculate extended energy

Calculate extended energy $E_{SQA}(u_s^*)$ with $x_s$

Calculate energy difference: $\Delta E_{SQA}(u_s^*) = E_{SQA}(u_s^*)^{new} - E_{SQA}(u_s^*)^{current}$

**If** $\Delta E_{SQA}(u_s^*) < 0$ or $\mathrm{rand}(0,1) < \exp\left(-\frac{\Delta E_{SQA}(u_s^*)}{T}\right)$

    Accept $u_s^*$ as current solution

**End if**

// Annealing schedules

Update temperature $T_{n+1} = vT_n, n = \{1, 2, ..., N_{iter}\}$, update tunnelling strength $\Gamma(n) = \Gamma_0 \cdot e^{-vn}$

// Convergence check

**If** maximum iterations $N_{iter}$ reached or objective function stabilises, then

    **Stop**

**End if**

**End for**

**Return** best solution: road link flows $x_s$, route flows $f_r^{\omega}$, $TSTT$, and optimal decision $u_s^*$.

Solving larger-scale networks using Quantum Annealing (QA) on D-Wave quantum hardware.

## Appendix C: Parameter value settings

**Table C1.** Road link parameters of the Nguyen-Dupuis network.

| Link No. | Connecting nodes | Road distance (km) | Free flow travel time (min) | Maximum speed (km/h) | Passing capacity (pcu/h) |
|---|---|---|---|---|---|
| 1 | 1-12 | 10.0 | 12 | 50 | 800 |
| 2 | 1-5 | 11.0 | 10 | 60 | 700 |
| 3 | 12-6 | 10.0 | 10 | 70 | 600 |
| 4 | 12-8 | 16.0 | 15 | 60 | 900 |
| 5 | 4-5 | 11.0 | 8 | 80 | 700 |
| 6 | 5-6 | 12.0 | 12 | 60 | 500 |
| 7 | 6-7 | 8.0 | 6 | 80 | 300 |
| 8 | 7-8 | 10.0 | 8 | 70 | 400 |
| 9 | 4-9 | 14.0 | 12 | 70 | 600 |
| 10 | 5-9 | 10.0 | 12 | 60 | 600 |
| 11 | 6-10 | 11.0 | 12 | 50 | 700 |
| 12 | 7-11 | 12.0 | 12 | 60 | 800 |
| 13 | 8-2 | 13.0 | 13 | 70 | 800 |
| 14 | 9-10 | 10.0 | 10 | 60 | 400 |
| 15 | 10-11 | 12.0 | 9 | 80 | 600 |
| 16 | 11-2 | 11.0 | 8 | 80 | 500 |
| 17 | 9-13 | 14.0 | 12 | 70 | 900 |
| 18 | 11-3 | 10.0 | 12 | 60 | 600 |
| 19 | 13-3 | 10.0 | 12 | 60 | 600 |

**Table C2.** OD demands and parameters.

| OD pair | Low-level demand (pcu/h) | Medium-level demand (pcu/h) | High-level demand (pcu/h) |
|---|---|---|---|
| 1 → 2 | 500 | 1000 | 1500 |
| 1 → 3 | 500 | 1000 | 1500 |
| 4 → 2 | 500 | 1000 | 1500 |
| 4 → 3 | 500 | 1000 | 1500 |

**Table C3.** Model and algorithm parameter settings.

| Parameter | Value (this study) | Parameter | Value (this study) |
|---|---|---|---|
| $\alpha$ | 0.15 | $c_s$ | Calculated via UE model (values reported in results section) |
| $\beta$ | 4.0 | $\beta_{st}$ | Calculated via UE model (values reported in results section) |
| $e$ | Randomly assigned to each disrupted road link from the interval [0.3, 0.7] | $T_0$ | 10 |
| $k$ | 2,3,4,5 | $\Gamma_0$ | 10 |
| $\lambda$ | 5000 ( $10 \times \max(c_s)$ ) | $v$ | 0.95 |
| $M$ | 10 | $N_{iter}$ | 50 |

For the ND test network, the maximum number of iterations was set to 50, which is sufficient to ensure convergence while avoiding unnecessary computation.

**Appendix D: Link impact and interaction coefficients in the ND network**

**Table D1.** Single-link impact coefficients under medium OD demand in ND network.

| Lin k | Remaining capacity ratio $e \in [0.3, 0.7]$ | Baseline $TSTT$ | Disrupted $TSTT$ | Single-link impact coefficient $c_s$ (Results rounded to two decimal places) |
|---|---|---|---|---|
| 1 | 0.578587674 | 5749.262154 | 6827.465377 | 1078.20 |
| 2 | 0.414455734 | 5749.262154 | 5873.286627 | 124.02 |
| 3 | 0.390740581 | 5749.262154 | 5812.972848 | 63.71 |
| 4 | 0.520525908 | 5749.262154 | 6272.085755 | 522.82 |
| 5 | 0.587787588 | 5749.262154 | 6178.627015 | 429.36 |
| 6 | 0.469242584 | 5749.262154 | 6461.474614 | 712.21 |
| 7 | 0.692305679 | 5749.262154 | 6429.222192 | 679.96 |
| 8 | 0.573931895 | 5749.262154 | 5749.480751 | 0.22 |
| 9 | 0.492372761 | 5749.262154 | 9506.183726 | 3756.92 |
| 10 | 0.456847007 | 5749.262154 | 6341.122796 | 591.86 |
| 11 | 0.437271206 | 5749.262154 | 5865.857227 | 116.60 |
| 12 | 0.591619883 | 5749.262154 | 5753.064254 | 3.802 |
| 13 | 0.475428898 | 5749.262154 | 6971.199272 | 1221.94 |
| 14 | 0.323871159 | 5749.262154 | 9290.127078 | 3540.86 |
| 15 | 0.459217702 | 5749.262154 | 7538.3722 | 1789.11 |
| 16 | 0.595198162 | 5749.262154 | 6273.770385 | 524.51 |
| 17 | 0.372996692 | 5749.262154 | 10775.24449 | 5025.98 |
| 18 | 0.370180702 | 5749.262154 | 7919.45321 | 2170.19 |
| 19 | 0.51262055 | 5749.262154 | 12428.65398 | 6679.39 |

**Table D2.** Interaction coefficients under medium OD demand in ND network.

| Roa d $s$ | Lin k $t$ | Baseline $TSTT$ | $TSTT_s$ (disrupted) | $TSTT_t$ (disrupted) | $TSTT_{st}$ (disrupted) | Interaction coefficient $\beta_{st}$ (Results rounded to two decimal places) |
|---|---|---|---|---|---|---|
| 1 | 2 | 5749.262 | 6827.465 | 5873.287 | 8609.724 | 1658.23 |
| 1 | 3 | 5749.262 | 6827.465 | 5812.973 | 6598.425 | -292.75 |
| 1 | 4 | 5749.262 | 6827.465 | 6272.086 | 7547.417 | 197.13 |
| 1 | 5 | 5749.262 | 6827.465 | 6178.627 | 8915.092 | 1658.26 |
| 1 | 6 | 5749.262 | 6827.465 | 6461.475 | 12553.98 | 5014.3 |
| 1 | 7 | 5749.262 | 6827.465 | 6429.222 | 18776.02 | 11268.59 |
| 1 | 8 | 5749.262 | 6827.465 | 5749.481 | 7369.393 | 541.71 |
| 1 | 9 | 5749.262 | 6827.465 | 9506.184 | 18931.11 | 8346.73 |
| 1 | 10 | 5749.262 | 6827.465 | 6341.123 | 9234.407 | 1815.08 |
| 1 | 11 | 5749.262 | 6827.465 | 5865.857 | 6558.613 | -385.45 |
| 1 | 12 | 5749.262 | 6827.465 | 5753.064 | 7743.332 | 912.06 |
| 1 | 13 | 5749.262 | 6827.465 | 6971.199 | 10736.43 | 2687.03 |
| 1 | 14 | 5749.262 | 6827.465 | 9290.127 | 8721.576 | -1646.75 |
| 1 | 15 | 5749.262 | 6827.465 | 7538.372 | 10339.71 | 1723.14 |
| 1 | 16 | 5749.262 | 6827.465 | 6273.77 | 7107.041 | -244.93 |
| 1 | 17 | 5749.262 | 6827.465 | 10775.24 | 16053.39 | 4199.94 |
| 1 | 18 | 5749.262 | 6827.465 | 7919.453 | 7090.772 | -1906.88 |
| 1 | 19 | 5749.262 | 6827.465 | 12428.65 | 24730.57 | 11223.72 |
| 2 | 3 | 5749.262 | 5873.287 | 5812.973 | 5981.822 | 44.82 |
| 2 | 4 | 5749.262 | 5873.287 | 6272.086 | 6078.181 | -317.93 |
| 2 | 5 | 5749.262 | 5873.287 | 6178.627 | 6225.535 | -77.12 |
| 2 | 6 | 5749.262 | 5873.287 | 6461.475 | 6227.723 | -357.78 |
| 2 | 7 | 5749.262 | 5873.287 | 6429.222 | 6884.035 | 330.79 |
| 2 | 8 | 5749.262 | 5873.287 | 5749.481 | 6136.536 | 263.03 |
| 2 | 9 | 5749.262 | 5873.287 | 9506.184 | 18262.25 | 8632.05 |
| 2 | 10 | 5749.262 | 5873.287 | 6341.123 | 7011.713 | 546.57 |
| 2 | 11 | 5749.262 | 5873.287 | 5865.857 | 5765.875 | -224.01 |
| 2 | 12 | 5749.262 | 5873.287 | 5753.064 | 5869.127 | -7.96 |
| 2 | 13 | 5749.262 | 5873.287 | 6971.199 | 7217.909 | 122.69 |
| 2 | 14 | 5749.262 | 5873.287 | 9290.127 | 8964.39 | -449.76 |
| 2 | 15 | 5749.262 | 5873.287 | 7538.372 | 12677.12 | 5014.73 |
| 2 | 16 | 5749.262 | 5873.287 | 6273.77 | 6664.852 | 267.06 |
| 2 | 17 | 5749.262 | 5873.287 | 10775.24 | 6297.992 | -4601.28 |
| 2 | 18 | 5749.262 | 5873.287 | 7919.453 | 7164.678 | -878.8 |
| 2 | 19 | 5749.262 | 5873.287 | 12428.65 | 19702.75 | 7150.07 |
| 3 | 4 | 5749.262 | 5812.973 | 6272.086 | 6143.849 | -191.95 |
| 3 | 5 | 5749.262 | 5812.973 | 6178.627 | 6844.392 | 602.05 |
| 3 | 6 | 5749.262 | 5812.973 | 6461.475 | 6242.519 | -282.67 |
| 3 | 7 | 5749.262 | 5812.973 | 6429.222 | 15990.36 | 9497.42 |
| 3 | 8 | 5749.262 | 5812.973 | 5749.481 | 5792.714 | -20.48 |
| 3 | 9 | 5749.262 | 5812.973 | 9506.184 | 9185.759 | -384.14 |
| 3 | 10 | 5749.262 | 5812.973 | 6341.123 | 6845.231 | 440.4 |
| 3 | 11 | 5749.262 | 5812.973 | 5865.857 | 5927.751 | -1.82 |
| 3 | 12 | 5749.262 | 5812.973 | 5753.064 | 5765.566 | -51.21 |
| 3 | 13 | 5749.262 | 5812.973 | 6971.199 | 6274.48 | -760.43 |
| 3 | 14 | 5749.262 | 5812.973 | 9290.127 | 7026.356 | -2327.48 |
| 3 | 15 | 5749.262 | 5812.973 | 7538.372 | 8584.875 | 982.79 |
| 3 | 16 | 5749.262 | 5812.973 | 6273.77 | 7042.718 | 705.24 |
| 3 | 17 | 5749.262 | 5812.973 | 10775.24 | 7106.387 | -3732.57 |
| 3 | 18 | 5749.262 | 5812.973 | 7919.453 | 7820.215 | -162.95 |
| 3 | 19 | 5749.262 | 5812.973 | 12428.65 | 8014.866 | -4477.5 |
| 4 | 5 | 5749.262 | 6272.086 | 6178.627 | 7818.301 | 1116.85 |
| 4 | 6 | 5749.262 | 6272.086 | 6461.475 | 6645.451 | -338.85 |
| 4 | 7 | 5749.262 | 6272.086 | 6429.222 | 11103.19 | 4151.14 |
| 4 | 8 | 5749.262 | 6272.086 | 5749.481 | 6252.348 | -19.96 |
| 4 | 9 | 5749.262 | 6272.086 | 9506.184 | 7634.364 | -2394.64 |
| 4 | 10 | 5749.262 | 6272.086 | 6341.123 | 7546.314 | 682.37 |
| 4 | 11 | 5749.262 | 6272.086 | 5865.857 | 6015.087 | -373.59 |
| 4 | 12 | 5749.262 | 6272.086 | 5753.064 | 7168.936 | 893.05 |
| 4 | 13 | 5749.262 | 6272.086 | 6971.199 | 8294.04 | 800.02 |
| 4 | 14 | 5749.262 | 6272.086 | 9290.127 | 7017.75 | -2795.2 |
| 4 | 15 | 5749.262 | 6272.086 | 7538.372 | 8012.716 | -48.48 |
| 4 | 16 | 5749.262 | 6272.086 | 6273.77 | 7072.104 | 275.51 |
| 4 | 17 | 5749.262 | 6272.086 | 10775.24 | 7062.15 | -4235.92 |
| 4 | 18 | 5749.262 | 6272.086 | 7919.453 | 7549.439 | -892.84 |
| 4 | 19 | 5749.262 | 6272.086 | 12428.65 | 33888.33 | 20936.86 |
| 5 | 6 | 5749.262 | 6178.627 | 6461.475 | 6939.392 | 48.55 |
| 5 | 7 | 5749.262 | 6178.627 | 6429.222 | 10416.92 | 3558.33 |
| 5 | 8 | 5749.262 | 6178.627 | 5749.481 | 7260.04 | 1081.19 |
| 5 | 9 | 5749.262 | 6178.627 | 9506.184 | 8007.691 | -1927.86 |
| 5 | 10 | 5749.262 | 6178.627 | 6341.123 | 6685.765 | -84.72 |
| 5 | 11 | 5749.262 | 6178.627 | 5865.857 | 6188.551 | -106.67 |
| 5 | 12 | 5749.262 | 6178.627 | 5753.064 | 7400.943 | 1218.51 |
| 5 | 13 | 5749.262 | 6178.627 | 6971.199 | 7684.781 | 284.22 |
| 5 | 14 | 5749.262 | 6178.627 | 9290.127 | 6356.294 | -3363.2 |
| 5 | 15 | 5749.262 | 6178.627 | 7538.372 | 7881.9 | -85.84 |
| 5 | 16 | 5749.262 | 6178.627 | 6273.77 | 8073.51 | 1370.37 |
| 5 | 17 | 5749.262 | 6178.627 | 10775.24 | 6534.999 | -4669.61 |
| 5 | 18 | 5749.262 | 6178.627 | 7919.453 | 8990.857 | 642.04 |
| 5 | 19 | 5749.262 | 6178.627 | 12428.65 | 20339.41 | 7481.39 |
| 6 | 7 | 5749.262 | 6461.475 | 6429.222 | 8707.772 | 1566.34 |
| 6 | 8 | 5749.262 | 6461.475 | 5749.481 | 7114.734 | 653.04 |
| 6 | 9 | 5749.262 | 6461.475 | 9506.184 | 11376.72 | 1158.33 |
| 6 | 10 | 5749.262 | 6461.475 | 6341.123 | 8781.72 | 1728.38 |
| 6 | 11 | 5749.262 | 6461.475 | 5865.857 | 6264.648 | -313.42 |
| 6 | 12 | 5749.262 | 6461.475 | 5753.064 | 6267.777 | -197.5 |
| 6 | 13 | 5749.262 | 6461.475 | 6971.199 | 7646.943 | -36.47 |
| 6 | 14 | 5749.262 | 6461.475 | 9290.127 | 11706.34 | 1704 |
| 6 | 15 | 5749.262 | 6461.475 | 7538.372 | 10013.27 | 1762.69 |
| 6 | 16 | 5749.262 | 6461.475 | 6273.77 | 7105.857 | 119.87 |
| 6 | 17 | 5749.262 | 6461.475 | 10775.24 | 9065.564 | -2421.89 |
| 6 | 18 | 5749.262 | 6461.475 | 7919.453 | 8061.223 | -570.44 |
| 6 | 19 | 5749.262 | 6461.475 | 12428.65 | 9029.374 | -4111.49 |
| 7 | 8 | 5749.262 | 6429.222 | 5749.481 | 7907.778 | 1478.34 |
| 7 | 9 | 5749.262 | 6429.222 | 9506.184 | 9881.492 | -304.65 |
| 7 | 10 | 5749.262 | 6429.222 | 6341.123 | 27941.28 | 20920.2 |
| 7 | 11 | 5749.262 | 6429.222 | 5865.857 | 8094.382 | 1548.56 |
| 7 | 12 | 5749.262 | 6429.222 | 5753.064 | 6904.655 | 471.63 |
| 7 | 13 | 5749.262 | 6429.222 | 6971.199 | 9136.922 | 1485.76 |
| 7 | 14 | 5749.262 | 6429.222 | 9290.127 | 12693.22 | 2723.14 |
| 7 | 15 | 5749.262 | 6429.222 | 7538.372 | 23179.43 | 14961.1 |
| 7 | 16 | 5749.262 | 6429.222 | 6273.77 | 9949.079 | 2995.35 |
| 7 | 17 | 5749.262 | 6429.222 | 10775.24 | 7091.561 | -4363.64 |
| 7 | 18 | 5749.262 | 6429.222 | 7919.453 | 32013.33 | 23413.91 |
| 7 | 19 | 5749.262 | 6429.222 | 12428.65 | 25487.24 | 12378.63 |
| 8 | 9 | 5749.262 | 5749.481 | 9506.184 | 15132.34 | 5625.94 |
| 8 | 10 | 5749.262 | 5749.481 | 6341.123 | 6571.365 | 230.02 |
| 8 | 11 | 5749.262 | 5749.481 | 5865.857 | 5777.408 | -88.67 |
| 8 | 12 | 5749.262 | 5749.481 | 5753.064 | 5844.496 | 91.21 |
| 8 | 13 | 5749.262 | 5749.481 | 6971.199 | 6771.648 | -199.77 |
| 8 | 14 | 5749.262 | 5749.481 | 9290.127 | 10394.9 | 1104.56 |
| 8 | 15 | 5749.262 | 5749.481 | 7538.372 | 6425.889 | -1112.7 |
| 8 | 16 | 5749.262 | 5749.481 | 6273.77 | 8501.904 | 2227.91 |
| 8 | 17 | 5749.262 | 5749.481 | 10775.24 | 6195.31 | -4580.15 |
| 8 | 18 | 5749.262 | 5749.481 | 7919.453 | 6167.7 | -1751.97 |
| 8 | 19 | 5749.262 | 5749.481 | 12428.65 | 27589.19 | 15160.31 |
| 9 | 10 | 5749.262 | 9506.184 | 6341.123 | 11746.14 | 1648.09 |
| 9 | 11 | 5749.262 | 9506.184 | 5865.857 | 10400.04 | 777.26 |
| 9 | 12 | 5749.262 | 9506.184 | 5753.064 | 6623.754 | -2886.23 |
| 9 | 13 | 5749.262 | 9506.184 | 6971.199 | 13702.35 | 2974.23 |
| 9 | 14 | 5749.262 | 9506.184 | 9290.127 | 16660.77 | 3613.72 |
| 9 | 15 | 5749.262 | 9506.184 | 7538.372 | 32366.45 | 21071.16 |
| 9 | 16 | 5749.262 | 9506.184 | 6273.77 | 22185.18 | 12154.49 |
| 9 | 17 | 5749.262 | 9506.184 | 10775.24 | 33826.44 | 19294.28 |
| 9 | 18 | 5749.262 | 9506.184 | 7919.453 | 24641.03 | 12964.65 |
| 9 | 19 | 5749.262 | 9506.184 | 12428.65 | 26560.02 | 10374.44 |
| 10 | 11 | 5749.262 | 6341.123 | 5865.857 | 7046.572 | 588.85 |
| 10 | 12 | 5749.262 | 6341.123 | 5753.064 | 6388.635 | 43.71 |
| 10 | 13 | 5749.262 | 6341.123 | 6971.199 | 7915.121 | 352.06 |
| 10 | 14 | 5749.262 | 6341.123 | 9290.127 | 7003.989 | -2878 |
| 10 | 15 | 5749.262 | 6341.123 | 7538.372 | 8901.968 | 771.73 |
| 10 | 16 | 5749.262 | 6341.123 | 6273.77 | 10450.77 | 3585.14 |
| 10 | 17 | 5749.262 | 6341.123 | 10775.24 | 14451.57 | 3084.47 |
| 10 | 18 | 5749.262 | 6341.123 | 7919.453 | 6314.966 | -2196.35 |
| 10 | 19 | 5749.262 | 6341.123 | 12428.65 | 9858.618 | -3161.9 |
| 11 | 12 | 5749.262 | 5865.857 | 5753.064 | 5803.678 | -65.98 |
| 11 | 13 | 5749.262 | 5865.857 | 6971.199 | 7388.141 | 300.35 |
| 11 | 14 | 5749.262 | 5865.857 | 9290.127 | 8220.416 | -1186.31 |
| 11 | 15 | 5749.262 | 5865.857 | 7538.372 | 6936.356 | -718.61 |
| 11 | 16 | 5749.262 | 5865.857 | 6273.77 | 8084.396 | 1694.03 |
| 11 | 17 | 5749.262 | 5865.857 | 10775.24 | 17050.16 | 6158.32 |
| 11 | 18 | 5749.262 | 5865.857 | 7919.453 | 7795.826 | -240.22 |
| 11 | 19 | 5749.262 | 5865.857 | 12428.65 | 29086.28 | 16541.03 |
| 12 | 13 | 5749.262 | 5753.064 | 6971.199 | 7273.487 | 298.49 |
| 12 | 14 | 5749.262 | 5753.064 | 9290.127 | 6351.277 | -2942.65 |
| 12 | 15 | 5749.262 | 5753.064 | 7538.372 | 8482.679 | 940.5 |
| 12 | 16 | 5749.262 | 5753.064 | 6273.77 | 7197.191 | 919.62 |
| 12 | 17 | 5749.262 | 5753.064 | 10775.24 | 7709.086 | -3069.96 |
| 12 | 18 | 5749.262 | 5753.064 | 7919.453 | 6309.685 | -1613.57 |
| 12 | 19 | 5749.262 | 5753.064 | 12428.65 | 11933.21 | -499.25 |
| 13 | 14 | 5749.262 | 6971.199 | 9290.127 | 7803.824 | -2708.24 |
| 13 | 15 | 5749.262 | 6971.199 | 7538.372 | 8452.448 | -307.86 |
| 13 | 16 | 5749.262 | 6971.199 | 6273.77 | 24389.53 | 16893.83 |
| 13 | 17 | 5749.262 | 6971.199 | 10775.24 | 9392.867 | -2604.31 |
| 13 | 18 | 5749.262 | 6971.199 | 7919.453 | 8999.562 | -141.83 |
| 13 | 19 | 5749.262 | 6971.199 | 12428.65 | 9871.572 | -3779.02 |
| 14 | 15 | 5749.262 | 9290.127 | 7538.372 | 8145.12 | -2934.12 |
| 14 | 16 | 5749.262 | 9290.127 | 6273.77 | 7380.086 | -2434.55 |
| 14 | 17 | 5749.262 | 9290.127 | 10775.24 | 10468.57 | -3847.54 |
| 14 | 18 | 5749.262 | 9290.127 | 7919.453 | 7225.607 | -4234.71 |
| 14 | 19 | 5749.262 | 9290.127 | 12428.65 | 11150.08 | -4819.44 |
| 15 | 16 | 5749.262 | 7538.372 | 6273.77 | 7789.467 | -273.41 |
| 15 | 17 | 5749.262 | 7538.372 | 10775.24 | 12064.28 | -500.07 |
| 15 | 18 | 5749.262 | 7538.372 | 7919.453 | 8719.494 | -989.07 |
| 15 | 19 | 5749.262 | 7538.372 | 12428.65 | 14089.62 | -128.14 |
| 16 | 17 | 5749.262 | 6273.77 | 10775.24 | 6973.792 | -4325.96 |
| 16 | 18 | 5749.262 | 6273.77 | 7919.453 | 7025.552 | -1418.41 |
| 16 | 19 | 5749.262 | 6273.77 | 12428.65 | 42980.59 | 30027.43 |
| 17 | 18 | 5749.262 | 10775.24 | 7919.453 | 13429.99 | 484.55 |
| 17 | 19 | 5749.262 | 10775.24 | 12428.65 | 18290.93 | 836.3 |
| 18 | 19 | 5749.262 | 7919.453 | 12428.65 | 14188.94 | -409.91 |

## Appendix E: Validation results on the larger transport network (Chicago Sketch)

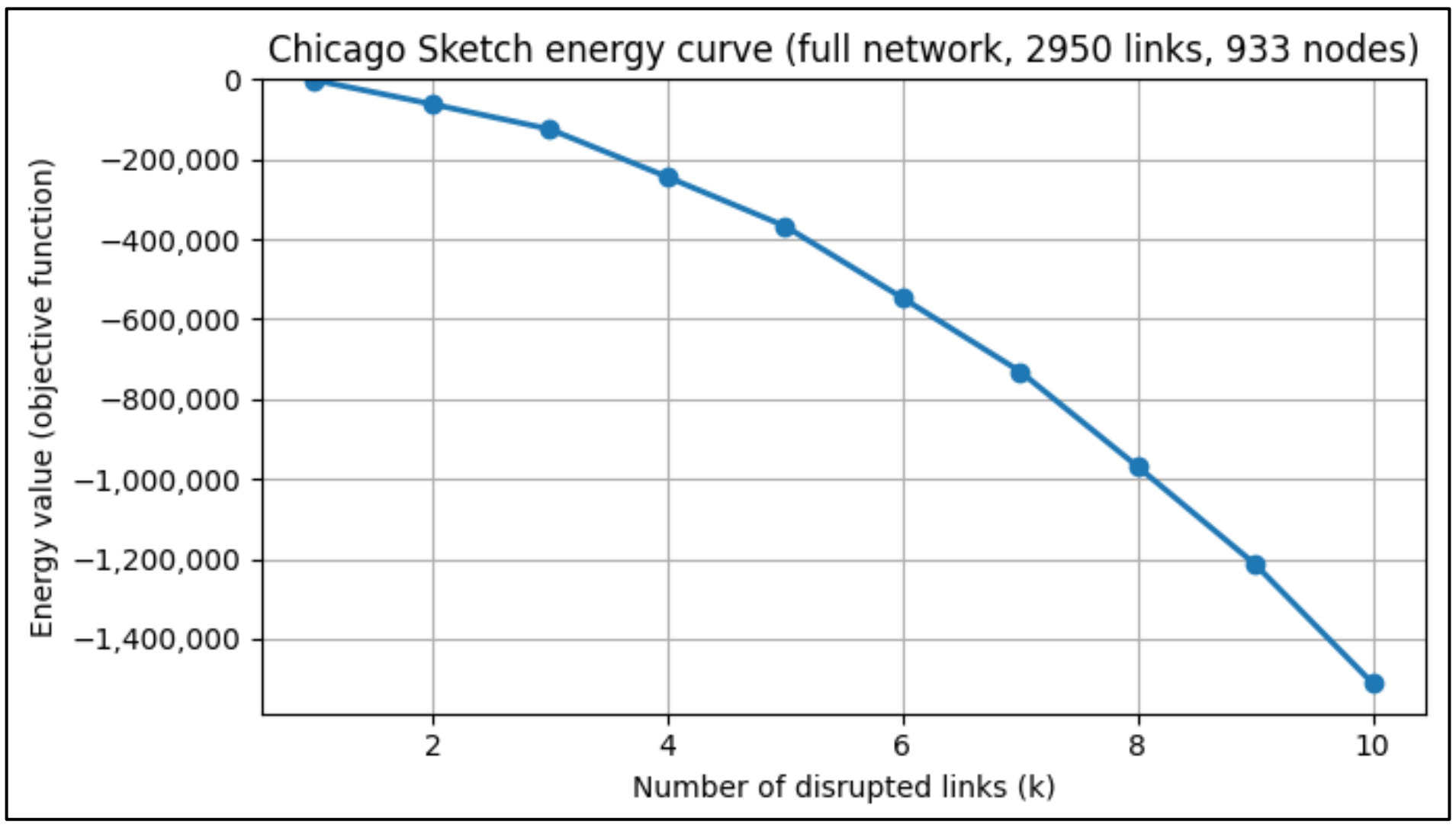

**Figure E1.** Energy curve for multiple disruptions in the Chicago Sketch network using D-Wave.

**Table E1.** D-Wave results for the Chicago Sketch network (2950 links).

| Number of disrupted links | Optimal energy | Runtime/s |
|---|---|---|
| 1 | -4000.946 | 70.3 |
| 2 | -62080.746 | 57.4 |
| 3 | -125450.567 | 63.1 |
| 4 | -244110.789 | 43.6 |
| 5 | -367458.902 | 76.6 |
| 6 | -547851.624 | 57.8 |
| 7 | -731790.062 | 58.1 |
| 8 | -969253.678 | 58.5 |
| 9 | -1213243.560 | 59.6 |
| 10 | -1521421.891 | 44.4 |

This Appendix presents the validation results based on the Chicago Sketch network (933 nodes and 2950 links), which are used to examine the computational feasibility and stability of the proposed quantum-classical hybrid optimisation framework under complex, large-scale transport networks. As shown in Figure E1, the system energy decreases monotonically with the increase in the number of disrupted links ( $k$ ), indicating that network vulnerability intensifies progressively with cumulative disruptions. All experiments were conducted on the D-Wave hybrid quantum solver, with an average runtime of approximately 50-70 seconds per run and a total completion time of less than 10 minutes. The results demonstrate that the framework can consistently obtain convergent solutions in the larger network, with smooth and steadily declining energy curves, further confirming the scalability and hardware-level feasibility of quantum computing for complex transport network vulnerability identification.

## Appendix F: Validation results on the larger transport network (Berlin Full)

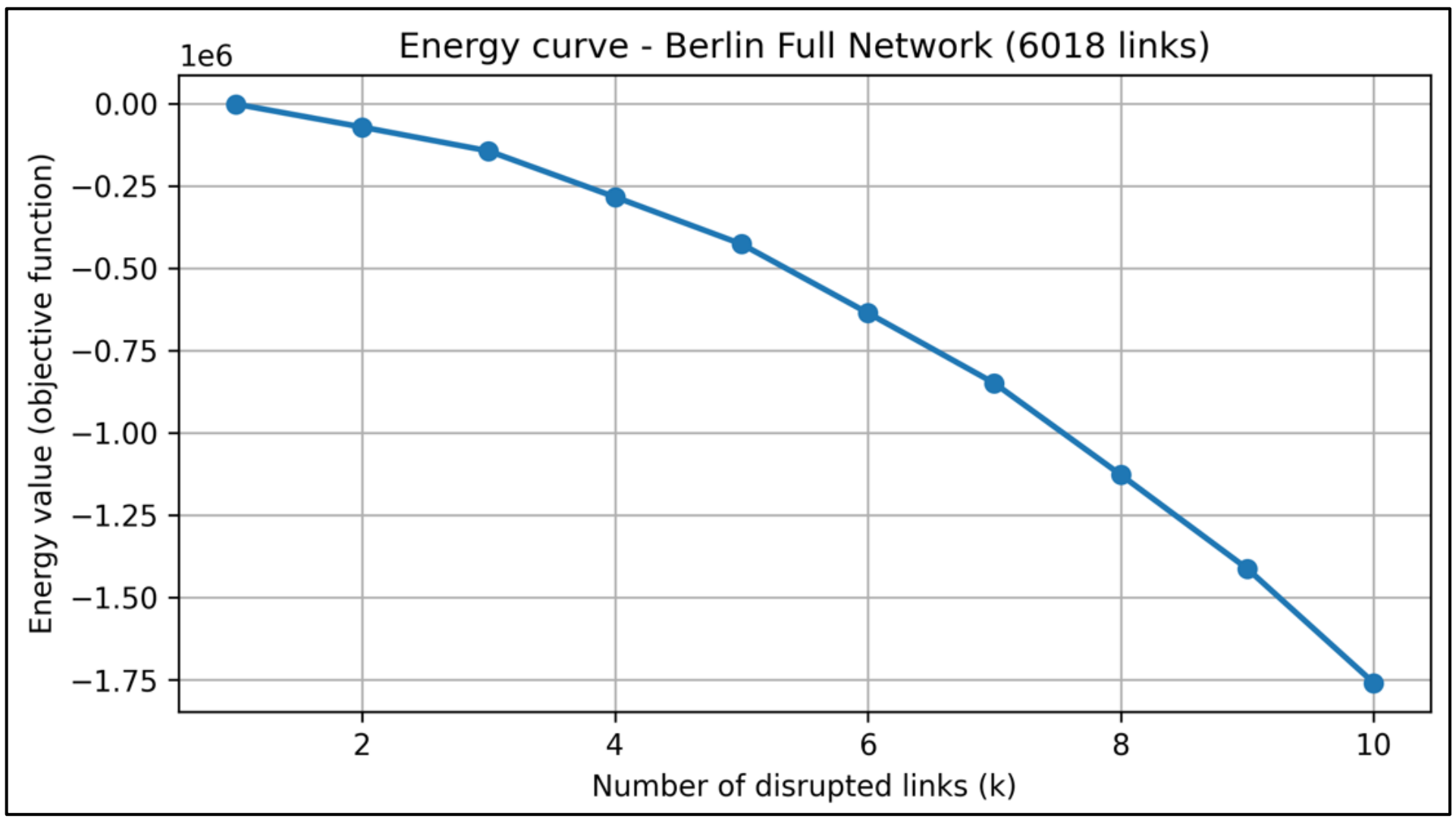

**Figure F1.** Energy curve for multiple disruptions in the Berlin Full network using D-Wave.

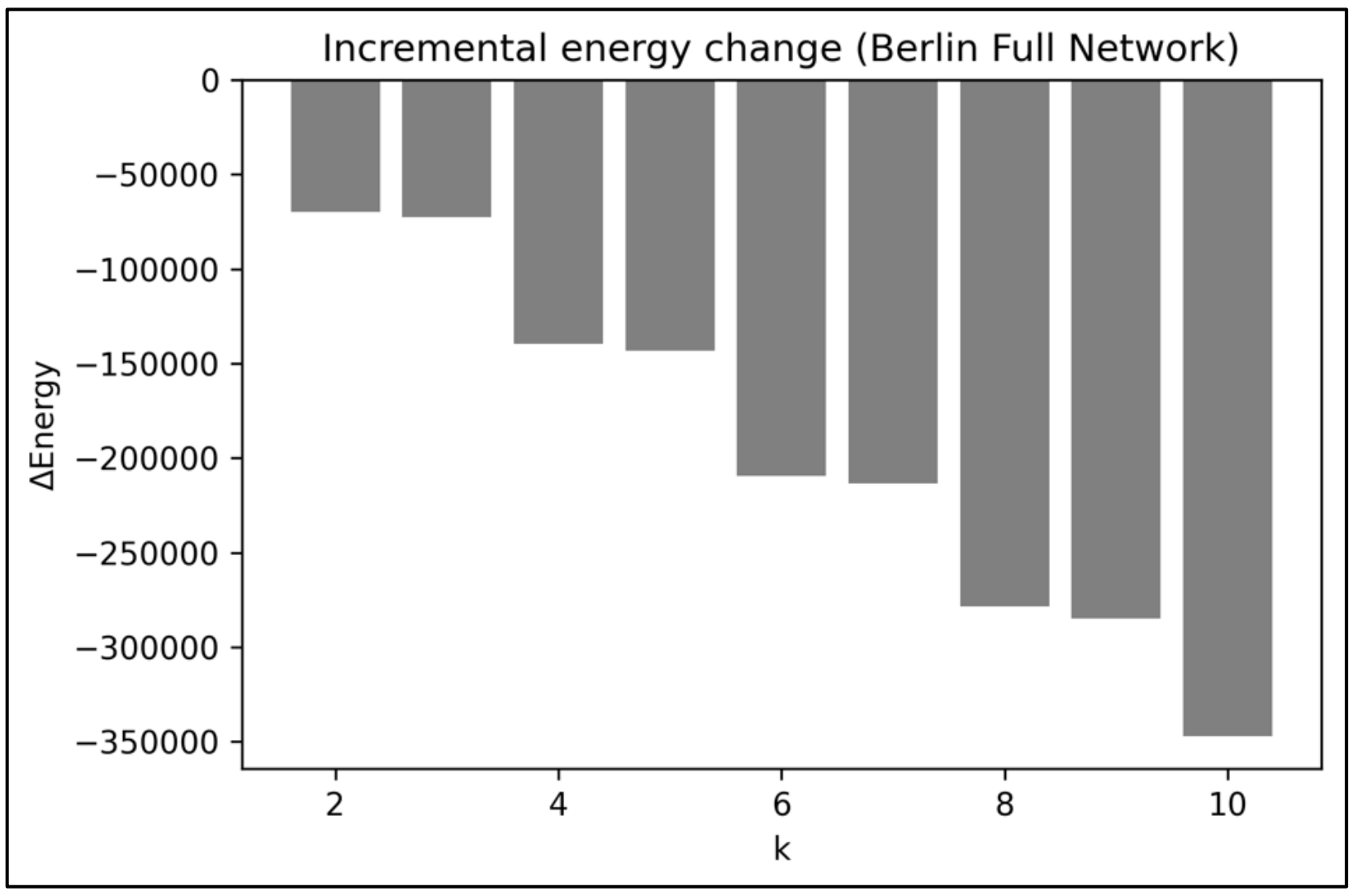

**Figure F2.** Incremental energy change in the Berlin Full network.

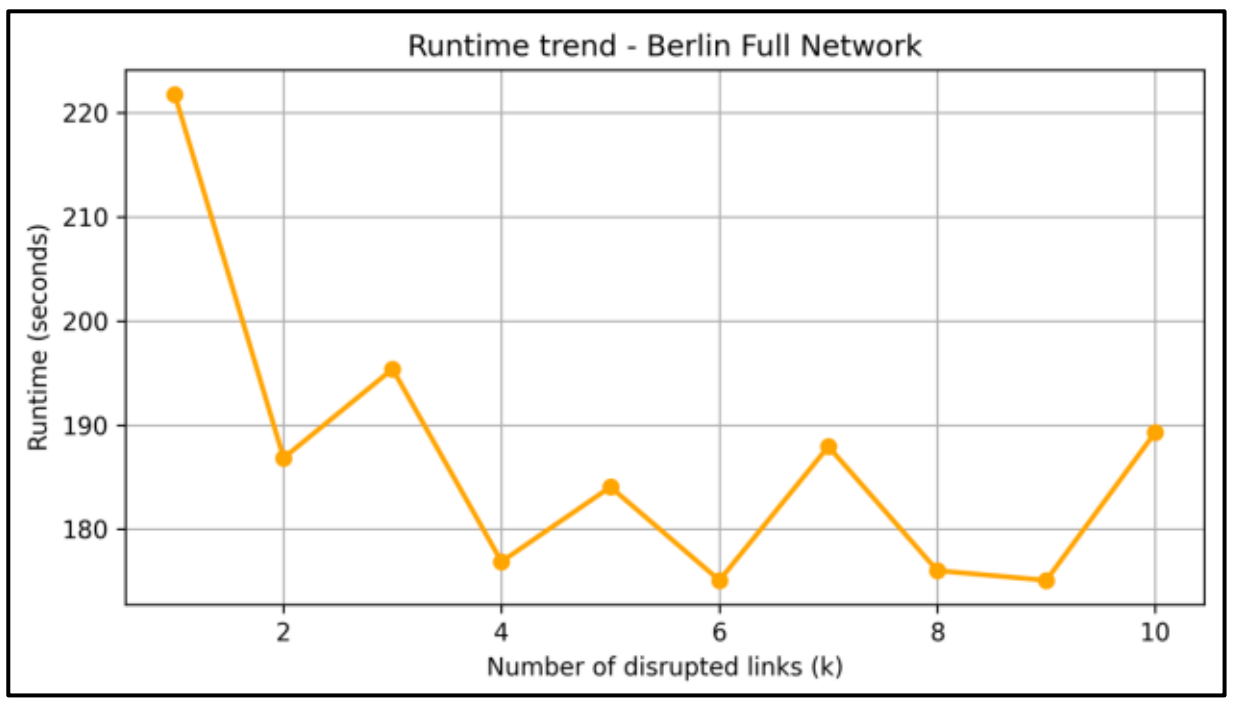

**Figure F3.** Runtime trend for Berlin Full Network.

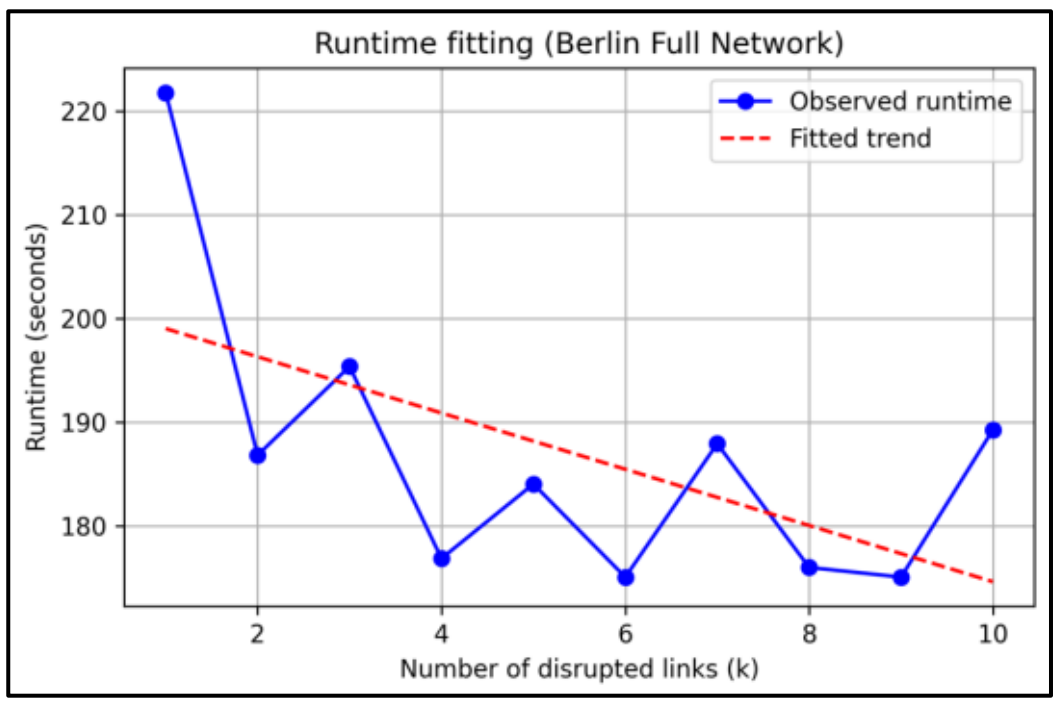

**Figure F4.** Runtime fitting for Berlin Full Network.

**Figure F1** presents the optimal energy variations of the Berlin full network under different disruption levels ($k=1-10$). The results show a continuous and nonlinear decline in system energy as the number of disrupted links increases, indicating that the network vulnerability intensifies cumulatively. In the early stage ($k \le 3$), the energy decreases gradually, suggesting that a few isolated link failures have limited impact. When the number of disruptions exceeds six, the energy drop accelerates, implying that the network enters a rapid degradation phase. Overall, the curve remains smooth and monotonically decreasing, demonstrating that the quantum-classical hybrid algorithm maintains good stability and convergence behaviour even in large-scale networks.

**Figure F2** depicts the incremental energy differences (ΔEnergy) between consecutive disruption levels. The results indicate that as $k$ increases, ΔEnergy becomes increasingly negative, meaning that each additional disrupted link produces a growing loss of system performance. The smaller ΔE values at the beginning imply that the network retains a certain redundancy capacity, while the sharp increase in ΔE after $k \geq 7$ marks the onset of high vulnerability. This pattern aligns with the well-known nonlinear fragility of transport networks: the system can absorb a limited number of disruptions, but once a critical threshold is exceeded, overall performance deteriorates rapidly.

**Figures F3** and **F4** illustrate the runtime behaviour of the D-Wave hybrid solver across different disruption levels and its fitted trend. The observed runtime remains within a narrow range of approximately 175 to 225 seconds, showing no significant increase as the number of disrupted links ( $k$ ) grows. This stability indicates that the hybrid algorithm scales efficiently, maintaining computational tractability even in large-scale QUBO optimisation. Minor fluctuations, such as those at $k=1$ and $k=3$ , are attributed to the solver's internal adaptive partitioning and parallel scheduling processes. The fitted trend line exhibits a nearly zero slope, confirming that the runtime varies linearly or remains almost constant with respect to $k$ . These results demonstrate that the hybrid solver effectively manages computational complexity through parallel decomposition and dynamic resource allocation, ensuring strong scalability and stable performance in the Berlin full network scenario.

**Table F1.** D-Wave results for the Berlin Full Network (6018 links).

| Number of disrupted links | Optimal energy | Runtime/s |
|---|---|---|
| 1 | -1998.945 | 221.8 |
| 2 | -71991.759 | 186.8 |
| 3 | -144524.886 | 195.4 |
| 4 | -284188.281 | 176.9 |
| 5 | -427209.721 | 184.1 |
| 6 | -636563.815 | 175.7 |
| 7 | -849883.696 | 187.9 |
| 8 | -1128236.951 | 176.0 |
| 9 | -1413179.483 | 175.1 |
| 10 | -1760121.540 | 189.3 |

**Appendix G: Parameter settings for metaheuristic algorithms used in the comparative experiments**

**Table G1.** Detailed table of parameter values and descriptions for metaheuristic algorithms.

| Algorithm | Parameter | Value / Description |
|---|---|---|
| GA | Population size | 200 |
| | Maximum generations | 300 |
| | Encoding scheme | Binary representation |
| | Selection method | Tournament selection |
| | Crossover probability | 0.50 |
| | Crossover operator | Two-point crossover |
| | Mutation probability | 0.50 |
| | Termination rule | Maximum generation budget reached |
| SA | Initial temperature | **10.0** (calibrated so that approximately 50 percent of uphill moves are accepted at the first iteration) |
| | Cooling schedule | Geometric cooling |
| | Cooling factor | 0.95 |
| | Temperature levels | 50 cooling stages |
| | Restart mechanism | Triggered by stagnation |
| | Termination rule | Maximum temperature cycle completed |
| PSO | Number of particles | 200 |
| | Maximum iterations | 300 |
| | Inertia weight | 0.7 |
| | Cognitive coefficient | 1.5 |
| | Social coefficient | 1.5 |
| | Position update rule | Binary mapping with sigmoid activation |
| | Velocity clamping | Enabled (max = 4.0) |
| | Termination rule | Maximum iteration budget reached |
| TS | Neighbourhood size | 100 candidate moves per iteration |
| | Move operator | Single link addition or removal |
| | Tabu tenure | 10 to 20 iterations (adaptive based on network size) |
| | Aspiration criterion | Accept tabu move if it yields improvement over the best-known solution |
| | Maximum iterations | 300 |
| | Maximum non improving iterations | 30 |
| | Restart mechanism | Restart triggered by stagnation |
| | Termination rule | Maximum iteration budget reached |

All parameter values are calibrated through pilot testing on the ND network to ensure that each heuristic algorithm achieves an average optimality gap below approximately 3% before being applied to larger networks (SF and Anaheim networks).

**Appendix H: Comparative runtime results of QA and heuristic algorithms across benchmark networks**

**Table H1.** Runtime results of QA and heuristic algorithms on the Sioux Falls network (Unit: s).

| $k$ | QA | GA | PSO | SA | TS |
|---|---|---|---|---|---|
| 1 | 4 | 5 | 6 | 4 | 6 |
| 2 | 3.9 | 7 | 8 | 5 | 7 |
| 3 | 4.2 | 9 | 11 | 6 | 9 |
| 4 | 4 | 12 | 14 | 8 | 12 |
| 5 | 4.3 | 15 | 18 | 10 | 16 |
| 6 | 3.9 | 19 | 23 | 13 | 21 |
| 7 | 3.9 | 24 | 29 | 17 | 27 |
| 8 | 3.9 | 30 | 36 | 22 | 34 |
| 9 | 3.8 | 37 | 45 | 28 | 43 |
| 10 | 3.9 | 45 | 55 | 35 | 54 |
| Total | 39.7 | 203 | 245 | 148 | 229 |

**Table H2.** Runtime results of QA and heuristic algorithms on the Anaheim network (Unit: min).

| $k$ | QA | GA | PSO | SA | TS |
|---|---|---|---|---|---|
| 1 | 0.27 | 6 | 7 | 5 | 6 |
| 2 | 0.27 | 8 | 10 | 7 | 9 |
| 3 | 0.28 | 11 | 14 | 10 | 13 |
| 4 | 0.29 | 15 | 19 | 14 | 18 |
| 5 | 0.28 | 20 | 25 | 19 | 24 |
| 6 | 0.33 | 26 | 32 | 25 | 31 |
| 7 | 0.28 | 33 | 40 | 32 | 39 |
| 8 | 0.27 | 45 | 54 | 43 | 52 |
| 9 | 0.28 | 65 | 80 | 60 | 69 |
| 10 | 0.27 | 85 | 110 | 72 | 90 |
| Total | 2.82 | 333 | 411 | 307 | 341 |

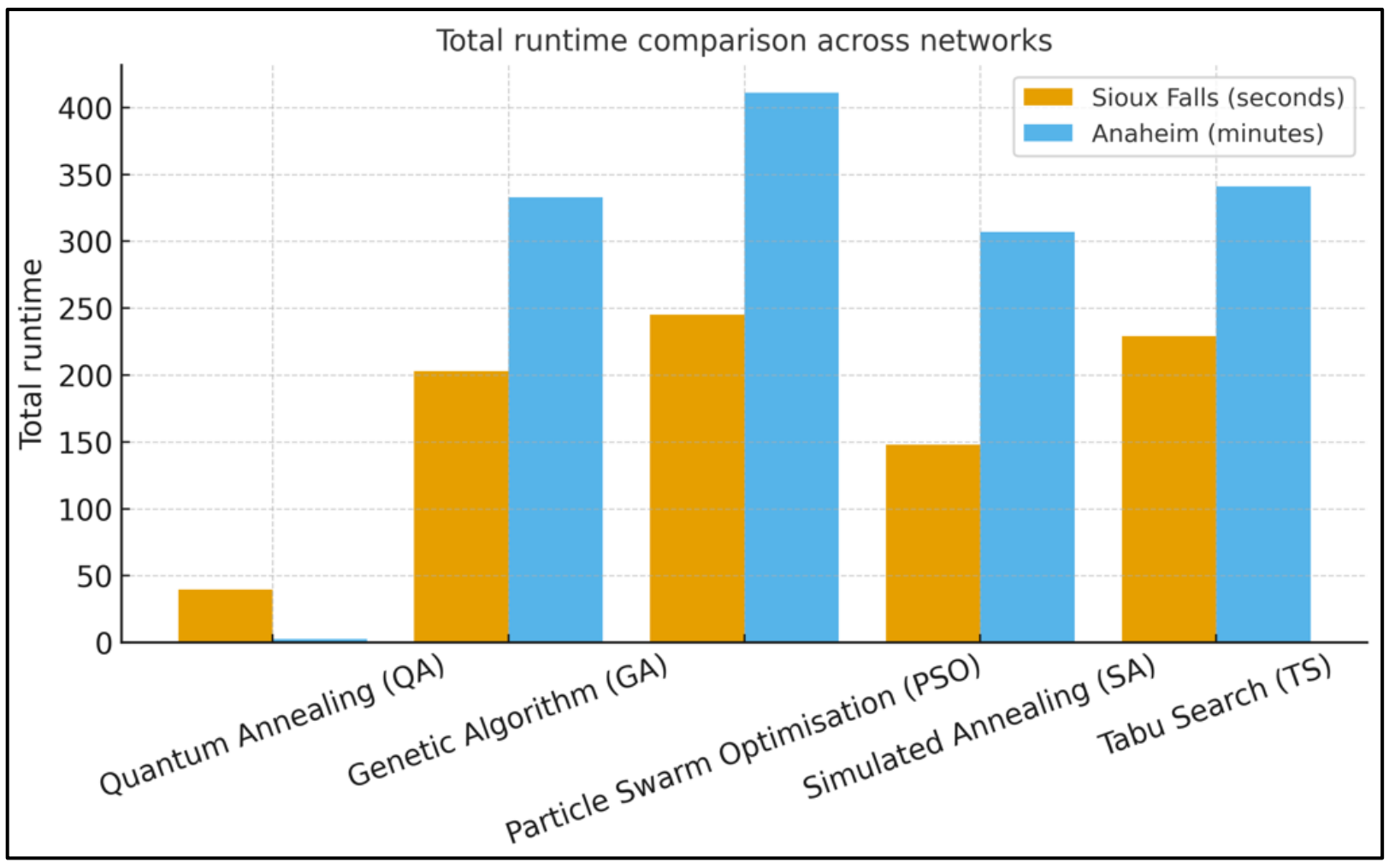

**Figure H1.** Total runtime comparison of QA and heuristic algorithms across SF and Anaheim networks.